\documentclass[12pt]{amsart}

\usepackage{amscd,latexsym,amsthm,amsfonts,amssymb,amsmath,amsxtra}
\usepackage[mathscr]{eucal}
\usepackage{hyperref}
\renewcommand\theequation{\thesection.\arabic{equation}}

\newcommand{\BA}{{\mathbb {A}}}

\newcommand{\BC}{{\mathbb {C}}}

\newcommand{\BR}{{\mathbb {R}}}

\newcommand{\CA}{{\mathcal {A}}}

\newcommand{\CC}{{\mathcal {C}}}
\newcommand{\CE}{{\mathcal {E}}}

\newcommand{\CH}{{\mathcal {H}}}

\newcommand{\CO}{{\mathcal {O}}}

\newcommand{\CS}{{\mathcal {S}}}

\newcommand{\CFJ}{{\mathcal {FJ}}}

\newcommand{\cusp}{{\mathrm{cusp}}}

\newcommand{\disc}{{\mathrm{disc}}}

\newcommand{\GL}{{\mathrm{GL}}}

\newcommand{\Ind}{{\mathrm{Ind}}}

\newcommand{\Mat}{{\mathrm{Mat}}}

\newcommand{\SO}{{\mathrm{SO}}}

\newcommand{\Sym}{{\mathrm{Sym}}}

\newcommand{\Sp}{{\mathrm{Sp}}}

\newcommand{\wt}{\widetilde}

\newcommand{\ol}{\overline}
\newcommand{\ul}{\underline}

\newcommand{\bs}{\backslash}

\def\bks{{\backslash}}

\def\diag{{\rm diag}}

\newtheorem{thm}{Theorem}[section]

\newtheorem{lem}[thm]{Lemma}

\newtheorem {conj}[thm]{Conjecture}

\newtheorem {ques/conj}[thm]{Question/Conjecture}

\newtheorem{defn}[thm]{Definition}

\makeatletter

\newcommand{\Rmnum}[1]{\expandafter\@slowromancap\romannumeral #1@}
\makeatother

\begin{document}
\renewcommand{\theequation}{\arabic{equation}}
\numberwithin{equation}{section}

\title[Fourier coefficients of residual representations]
{On Fourier coefficients of certain residual representations of symplectic groups}

\author{Dihua Jiang}
\address{School of Mathematics\\
University of Minnesota\\
Minneapolis, MN 55455, USA}
\email{dhjiang@math.umn.edu}

\author{Baiying Liu}
\address{Department of Mathematics\\
University of Utah\\
155 S 1400 E Room 233, Salt Lake City, UT 84112-0090, USA.}
\email{liu@math.utah.edu}


\date{\today}


\thanks{The research of the first named author is supported in part by the NSF Grants DMS--1301567, and the research of the second
named author is supported in part by NSF Grants DMS-1302122, and in part by a postdoc research fund from Department of Mathematics, University of Utah}

\begin{abstract}
In the theory of automorphic descents developed by Ginzburg, Rallis and Soudry in \cite{GRS11}, the structure of Fourier coefficients of
the residual representations of certain special Eisenstein series plays an essential role. Started from \cite{JLZ13}, the authors are looking for
more general residual representations, which may yield more general theory of automorphic descents. In this paper, we investigate the structure
of Fourier coefficients of certain residual representations of symplectic groups, associated with certain interesting families of global Arthur parameters. On one hand, the results partially confirm a conjecture proposed by the first named author in \cite{J14} on relations between the
global Arthur parameters and the structure of Fourier coefficients of the automorphic representations in the associated global Arthur packets.
On the other hand, the results of this paper can be regarded as a first step towards more general automorphic descents for
symplectic groups, which will be considered in our future work.
\end{abstract}

\maketitle



\section{Introduction}

Let $\Sp_{2n}$ be the symplectic group with symplectic form
$$\begin{pmatrix}
0 & v_n\\
-v_n & 0
\end{pmatrix},$$
where $v_n$ is an $n \times n$ matrix with $1$'s on the second diagonal and $0$'s elsewhere. Fix a Borel subgroup $B=TU$ of $\Sp_{2n}$,
where the maximal torus $T$ consists of elements of the following form
$$
\diag(t_1,\cdots,t_n;t_n^{-1},\cdots,t_1^{-1})
$$
and the unipotent radical $U$ consists of all upper unipotent matrices in $\Sp_{2n}$. Let $F$ be a number field and $\BA$ be the ring of
adeles of $F$.

The structure of Fourier coefficients for the residual representations of $\Sp_{4n}(\BA)$, with cuspidal support $(\GL_{2n},\tau)$, played an
indispensable role in the theory of automorphic descent from $\GL_{2n}$ to the metaplectic double cover of $\Sp_{2n}$ by Ginzburg, Rallis and
Soudry in \cite{GRS11}. As tested in a special case in our recent work joint with Xu and Zhang in \cite{JLXZ14}, we expected the residual
representations investigated in \cite{JLZ13} may play important roles in extending the theory of automorphic descent in \cite{GRS11} to a more
general setting. In this paper, we take certain interesting families of residual representations of $\Sp_{2n}(\BA)$ obtained in \cite{JLZ13} and
study the structure of their Fourier coefficients associated to nilpotent orbits as described in \cite{J14}.
On one hand, the results of this paper partially confirm a conjecture proposed by the first named author in \cite{J14} on relations between the
global Arthur parameters and the structure of Fourier coefficients of the automorphic representations in the corresponding global Arthur packets.
On the other hand, these results are preliminary steps towards the theory of more general automorphic descents for
symplectic groups, which will be considered in our future work.

We first recall the global Arthur parameters for $\Sp_{2n}$ and the discrete spectrum, and the conjecture made in \cite{J14}.
Then we recall what have been proved about this conjecture before this current paper, in particular
the results obtained in \cite{JL15c}. Finally we describe more explicitly the objective of this paper. The main results will be precisely
stated in Section 2.

\subsection{Arthur parameters and the discrete spectrum}
Let $F$ be a number field and $\BA$ be the ring of adeles of $F$.
Recall that the dual group of $G_n=\Sp_{2n}$ is $\SO_{2n+1}(\BC)$. The set of global Arthur parameters for the discrete spectrum of
the space of all square-integrable automorphic functions on $\Sp_{2n}(\BA)$ is denoted by $\wt{\Psi}_2(\Sp_{2n})$, following the notation in
\cite{Ar13}. The elements of $\wt{\Psi}_2(\Sp_{2n})$ are of the form
\begin{equation}\label{psi1}
\psi:=\psi_1\boxplus\psi_2\boxplus\cdots\boxplus\psi_r,
\end{equation}
where $\psi_i$ are pairwise distinct simple global Arthur parameters of orthogonal type. A simple global Arthur parameter is formally given
by $(\tau,b)$ with an integer $b\geq 1$, and with $\tau\in\CA_\cusp(a)$ being an irreducible unitary cuspidal automorphic representation of
$\GL_a(\BA)$.

In the expression \eqref{psi1}, one has that
$\psi_i=(\tau_i,b_i)$ with $\tau_i\in\CA_\cusp(a_i)$, $2n+1 = \sum_{i=1}^r a_ib_i$,
and $\prod_i \omega_{\tau_i}^{b_i} = 1$ (the condition
on the central character of the parameter), following \cite[Section 1.4]{Ar13}. In order for all $\psi_i$s to be of orthogonal type, the
simple parameters $\psi_i=(\tau_i,b_i)$ for $i=1,2,\cdots,r$
satisfy the following parity condition:
if $\tau_i$ is of symplectic type (i.e., $L(s, \tau_i, \wedge^2)$ has a pole at $s=1$), then $b_i$ is even; and
if $\tau_i$ is of orthogonal type (i.e., $L(s, \tau_i, \Sym^2)$ has a pole at $s=1$), then $b_i$ is odd.
A global Arthur parameter $\psi=\boxplus_{i=1}^r (\tau_i,b_i)$ is called generic if $b_i=1$ for all $1 \leq i \leq r$.

\begin{thm}[Theorem 1.5.2, \cite{Ar13}] For each global Arthur parameter $\psi\in\wt{\Psi}_2(\Sp_{2n})$, there exists a global Arthur
packet $\wt{\Pi}_\psi$. The discrete spectrum of $\Sp_{2n}(\BA)$ has the following decomposition
$$
L^2_\disc(\Sp_{2n}(F)\bks\Sp_{2n}(\BA))
\cong\oplus_{\psi\in\wt{\Psi}_2(\Sp_{2n})}\oplus_{\pi\in\wt{\Pi}_\psi(\epsilon_\psi)}\pi,
$$
where $\wt{\Pi}_\psi(\epsilon_\psi)$ denotes the subset of $\wt{\Pi}_\psi$ consisting of members which occur in the
discrete spectrum of $\Sp_{2n}(\BA)$.
\end{thm}

\subsection{A conjecture on the Fourier coefficients}
We will use the notation in \cite{JL15a} and \cite{JL15c} freely. Following \cite[Section2]{JL15a}, for a symplectic partition
$\ul{p}$ of $2n$, or equivalently each $F$-stable unipotent orbit $\CO_{\ul{p}}$,
via the standard $\mathfrak{sl}_2(F)$-triple, one may construct an $F$-unipotent subgroup $V_{\ul{p},2}$. In this case, the $F$-rational unipotent
orbits in the $F$-stable unipotent orbit $\CO_{\ul{p}}$ are parameterized by a data $\ul{a}$ (see \cite[Section 2]{JL15a} for detail),
which defines a character $\psi_{\ul{p},\ul{a}}$ of $V_{\ul{p},2}(\BA)$. This character $\psi_{\ul{p},\ul{a}}$ is automorphic in the sense that
it is trivial on $V_{\ul{p},2}(F)$.
The $\psi_{\underline{p}, \underline{a}}$-Fourier coefficient
of an automorphic form $\varphi$ on $\Sp_{2n}(\BA)$ is defined by
\begin{equation}\label{fc}
\varphi^{\psi_{\underline{p}, \underline{a}}}(g):=\int_{V_{\underline{p}, 2}(F)\bks V_{\underline{p}, 2}(\BA)}
\varphi(vg) \psi_{\underline{p}, \underline{a}}(v)^{-1} dv.
\end{equation}
We say that an irreducible automorphic representation $\pi$ of $\Sp_{2n}(\BA)$
has a nonzero $\psi_{\underline{p}, \underline{a}}$-Fourier coefficient or
a nonzero Fourier coefficient attached to a (symplectic) partition $\ul{p}$ if there exists an
automorphic form $\varphi$ in the space of $\pi$ with a nonzero $\psi_{\underline{p}, \underline{a}}$-Fourier coefficient
$\varphi^{\psi_{\underline{p}, \underline{a}}}(g)$, for some choice of $\ul{a}$.
For any irreducible automorphic representation $\pi$
of $\Sp_{2n}(\BA)$, as in \cite{J14}, we define $\frak{p}^m(\pi)$ (which corresponds to $\mathfrak{n}^m(\pi)$ in the notation of \cite{J14})
to be the set of all symplectic partitions $\underline{p}$ with the properties
that $\pi$ has a nonzero $\psi_{\underline{p}, \underline{a}}$-Fourier coefficient
for some choice of $\underline{a}$, but for any ${\underline{p}'} > \underline{p}$ (with
the natural ordering of partitions), $\pi$ has no nonzero Fourier coefficients
attached to ${\underline{p}'}$. It is generally believed (and may be called a conjecture) that the set $\frak{p}^m(\pi)$ contains
only one partition for any irreducible automorphic representation $\pi$ (or locally for any irreducible admissible representation $\pi$).
We refer to Section 3, Conjecture 3.1 in particular, of \cite{JL15b} for more detailed discussions on this issue.

As in \cite{J14}, $\wt{\Pi}_\psi(\epsilon_\psi)$ is called the automorphic $L^2$-packet attached to the global Arthur parameter $\psi$.
For each $\psi$ of the form in \eqref{psi1}, let $\underline{p}(\psi)=[(b_1)^{(a_1)}\cdots(b_r)^{(a_r)}]$ be a partition of $2n+1$ attached to the
global Arthur parameter $\psi$, following the discussion in \cite[Section 4]{J14}.
For $\pi\in\wt{\Pi}_\psi(\epsilon_\psi)$, the structure of the global Arthur parameter $\psi$ deduces constraints on
the structure of $\frak{p}^m(\pi)$, which are given by the following conjecture.

\begin{conj}[Conjecture 4.2, \cite{J14}]\label{cubmfc}
For any $\psi\in\wt{\Psi}_2(\Sp_{2n})$, let $\wt{\Pi}_{\psi}(\epsilon_\psi)$ be the automorphic $L^2$-packet attached to $\psi$.
Then the following hold.
\begin{enumerate}
\item[(1)] Any symplectic partition $\ul{p}$ of $2n$, if $\ul{p}>\eta_{{\frak{g}^\vee,\frak{g}}}(\underline{p}(\psi))$, does
not belong to $\frak{p}^m(\pi)$ for any $\pi\in\wt{\Pi}_{\psi}(\epsilon_\psi)$.
\item[(2)] For every $\pi\in\wt{\Pi}_{\psi}(\epsilon_\psi)$, any partition $\ul{p}\in\frak{p}^m(\pi)$ has the property that
$\ul{p}\leq \eta_{{\frak{g}^\vee,\frak{g}}}(\ul{p}(\psi))$.
\item[(3)] There exists at least one member $\pi\in\wt{\Pi}_{\psi}(\epsilon_\psi)$ having the property that
$\eta_{{\frak{g}^\vee,\frak{g}}}(\ul{p}(\psi))\in \frak{p}^m(\pi)$.
\end{enumerate}
Here $\eta_{{\frak{g}^\vee,\frak{g}}}$ denotes the Barbasch-Vogan duality map
$($see Definition \ref{def1}$)$ from the partitions for $\frak{so}_{2n+1}(\BC)$ to
the partitions for $\frak{sp}_{2n}(\BC)$.
\end{conj}

We remark that Part (2) is stronger than Part (1) in Conjecture \ref{cubmfc}. More related discussions can be found in \cite{JL15b}.

There have been progress toward the proof of Conjecture \ref{cubmfc}. When the global Arthur parameter $\psi=\boxplus_{i=1}^r (\tau_i,1)$ is
generic, in Conjecture \ref{cubmfc}, Part (1) is trivial and Part (2) is automatic; and Part (3) of Conjecture \ref{cubmfc} can be viewed as
the global version of the Shahidi conjecture,
namely, any global tempered L-packet has a generic member. This can be proved following the theory of automorphic descent developed by
Ginzburg, Rallis and Soudry (\cite{GRS11}) and the endoscopy classification of Arthur (\cite{Ar13}). We refer to Section 3.1, Theorem 3.3
in particular, of \cite{JL15b} for more precise discussion on this issue.
Hence Conjecture \ref{cubmfc} holds
for all generic global Arthur parameters, and those $\pi$ satisfying Part (3) are generic cuspidal representations.

For Arthur parameters of form $\psi=(\tau,b)\boxplus (1_{\GL_1(\BA)}, 1)$, where $\tau$ is an irreducible cuspidal representation of
$\GL_{2k}(\BA)$ and is of symplectic type, and $b$ is even, one has that $\ul{p}(\psi)=[b^{(2k)}1]$. In this case,
Part (3) of Conjecture \ref{cubmfc} has been proved by the second named author in \cite{L13b}, where it is also shown that
$\frak{p}^m(\pi)$ contains only one partition in this particular case.

For a general global Arthur parameter $\psi$, Part (1) of Conjecture \ref{cubmfc} is completely proved in \cite{JL15c}.
We remark that if we assume that $\frak{p}^m(\pi)$ contains only one partition, then Part (2) of Conjecture \ref{cubmfc}
essentially follows from Parts (1) and (3) of Conjecture \ref{cubmfc} plus certain local constraints at unramified local places as
discussed in \cite{JL15c}. We omit the details here. However, without knowing that the set $\frak{p}^m(\pi)$ contains only one partition,
Part (2) of Conjecture \ref{cubmfc} is also settled in \cite{JL15c} partially, namely, any symplectic partition $\ul{p}$ of $2n$, if $\ul{p}$ is bigger than $\eta_{{\frak{g}^\vee,\frak{g}}}(\underline{p}(\psi))$ under the lexicographical ordering, does
not belong to $\frak{p}^m(\pi)$ for any $\pi\in\wt{\Pi}_{\psi}(\epsilon_\psi)$.
We refer to \cite[Section 4]{J14} and also \cite{JL15b} for more discussion on this conjecture and related topics.

\subsection{The objective of this paper}
We start here to investigate Part (3) of Conjecture \ref{cubmfc}. This means that we have to construct or determine a particular
member in a given automorphic $L^2$-packet $\wt{\Pi}_\psi(\epsilon_\psi)$ attached to a general global Arthur parameter $\psi$, whose
Fourier coefficients achieve the partition $\eta_{{\frak{g}^\vee,\frak{g}}}(\ul{p}(\psi))$. We expect that such members should be the
distinguished members in $\wt{\Pi}_\psi(\epsilon_\psi)$, following the Whittaker normalization in the sense of Arthur for global generic
Arthur parameters (\cite{Ar13}). For general non-generic global Arthur parameters, the distinguished members in $\wt{\Pi}_\psi(\epsilon_\psi)$ can be
certain residual representations determined by $\psi$ as conjectured by M{\oe}glin in \cite{M08} and \cite{M11}, or
certain cuspidal automorphic representations, which may be explicitly constructed through the framework of endoscopy
correspondences as outlined in \cite{J14}. Due to the different nature of the two construction methods, we are going to treat them separately, in order to prove Part (3) of Conjecture \ref{cubmfc}.

As explained in \cite{JL15b}, when the distinguished members $\pi$ in a given $\wt{\Pi}_\psi(\epsilon_\psi)$ are residual representations,
they can be constructed explicitly from the given cuspidal data. In this case, our method is to establish
the nonvanishing of the Fourier coefficients of those $\pi$ associated to the partition $\eta_{{\frak{g}^\vee,\frak{g}}}(\ul{p}(\psi))$,
in terms of
the nonvanishing condition (Fourier coefficients or periods) on the construction data that is also defined by the given non-generic global
Arthur parameter $\psi$. Hence such a method can be regarded as a natural extension of the well-known
Langlands-Shahidi method from generic Eisenstein series (\cite{Sh10}) to non-generic Eisenstein series, and in  particular to the singularity
of Eisenstein series, i.e. the residues of Eisenstein series. On the other hand, this method can also be regarded as an extension of the
automorphic descent method of Ginzburg-Rallis-Soudry for a particular residual representations (\cite{GRS11}) to general residual representations.

In this paper, we are going to test our method for these non-generic global Arthur parameters $\psi$, whose
automorphic $L^2$-packets $\wt{\Pi}_\psi(\epsilon_\psi)$ contain the residual representations that are completely determined
in our previous work joint with Zhang (\cite{JLZ13}). Those non-generic global Arthur parameters of $\Sp_{2n}(\BA)$ are of the following form
$$
\psi=(\tau_1,b_1) \boxplus \boxplus_{i=2}^r (\tau_i,1),\ \text{with}\  b_1 > 1,
$$
which has three cases, depending on the symmetry of $\tau_1$ and the relation of $\tau_1$ with $\tau_i$ for $i=2,3,\cdots,r$.

\textbf{Case \Rmnum{1}:} $\psi=(\tau, 2b+1) \boxplus \boxplus_{i=2}^r (\tau_i,1),$
where $b \geq 1$, $\tau \ncong \tau_i$, for any $2 \leq i \leq r$.

\textbf{Case \Rmnum{2}:} $\psi=(\tau, 2b+1) \boxplus (\tau, 1)
\boxplus \boxplus_{i=3}^r (\tau_i,1),$
where $b \geq 1$, $\tau \ncong \tau_i$, for any $3 \leq i \leq r$.

\textbf{Case \Rmnum{3}:} $\psi=(\tau, 2b) \boxplus \boxplus_{i=2}^r (\tau_i,1),$
where $b \geq 1$.

For $\psi\in\wt{\Psi}_2(\Sp_{2n})$, $\tau\in\CA_\cusp(\GL_a)$
is of orthogonal type in \textbf{Case \Rmnum{1}} and \textbf{Case \Rmnum{2}}, and of symplectic type in \textbf{Case \Rmnum{3}}. Of course,
the remaining $\tau_i$s are of orthogonal type in all three cases.

When $\tau$ is of orthogonal type, i.e. in both \textbf{Case \Rmnum{1}} and \textbf{Case \Rmnum{2}}, the corresponding
residual representations given in \cite{JLZ13} must be nonzero. We prove in this paper
Part (3) of Conjecture \ref{cubmfc} those two cases, and refer to Section 2 for more details.

When $\tau$ is of symplectic type and $r\geq 2$, the relation between
$\tau$ and $\tau_i$, for $i=2,3,\cdots,r$, is governed by the corresponding Gan-Gross-Prasad conjecture (\cite{GGP12}),
which controls the structure of
the automorphic $L^2$-packet $\wt{\Pi}_{\psi}(\epsilon_\psi)$. We prove Part (3) of Conjecture \ref{cubmfc} for \textbf{Case \Rmnum{3}} when
$\wt{\Pi}_{\psi}(\epsilon_\psi)$ contains residual representations. While the automorphic $L^2$-packet $\wt{\Pi}_{\psi}(\epsilon_\psi)$
does not contain any residual representation, the situation is more involved, and will be left for a separate treatment in our future work.
We discuss with more details in Section 2.

We will state the main results more explicitly in Section 2. After recalling a technical lemma from \cite{JL15b} in Section 3, we are ready to treat
\textbf{Case \Rmnum{1}} in both Sections 4 and 5. \textbf{Case \Rmnum{2}} is treated in Section 6. The final section is devoted to
\textbf{Case \Rmnum{3}}. One may find more detailed description of the arguments and methods used in the proof of those cases in
each relevant section.

We would like to thank David Soudry for helpful discussion on related topics. We also would like to thank the referee for careful reading of the paper and helpful comments.

\section{The Main Results}

After introducing more notation and basic facts about discrete spectrum and Fourier coefficients attached to partitions, we will state
the main results explicitly for each case.

Throughout the paper, we let $P^{2n}_r = M^{2n}_r N^{2n}_r$ (with $1 \leq r \leq n$) be the standard parabolic subgroup
of ${\Sp}_{2n}$ with Levi part $M^{2n}_r$ isomorphic to $\GL_r \times \Sp_{2n-2r}$,
$N^{2n}_r$ is the unipotent radical. And let
$\wt{P}^{2n}_r(\BA)=\wt{M}^{2n}_r(\BA) N^{2n}_r(\BA)$ be the pre-image of $P^{2n}_r(\BA)$ in $\wt{\Sp}_{2n}(\BA)$ (the superscript $2n$ may be dropped when there is no confusion). The description of the three cases was briefly given in \cite{JL15b}. Here are the details.

\subsection{\textbf{Case \Rmnum{1}}}
 $\psi\in\wt{\Psi}_2(\Sp_{2n})$ is written as
 \begin{equation}\label{case1psi}
 \psi=(\tau, 2b+1) \boxplus \boxplus_{i=2}^r (\tau_i,1),
 \end{equation}
where $b \geq 1$, $\tau \ncong \tau_i$, for any $2 \leq i \leq r$. Assume $\tau\in\CA_\cusp(\GL_a)$ has central character $\omega_{\tau}$, and
$\tau_i\in\CA_\cusp(\GL_{a_i})$ has central character $\omega_{\tau_i}$, $2 \leq i \leq r$.
Following the definition of $\wt{\Psi}_2(\Sp_{2n})$, one must have that $2n+1=a(2b+1)+\sum_{i=2}^r {a_i}$,
and $\omega_{\tau}^{2b+1} \cdot \prod_{i=2}^r \omega_{\tau_i} = 1$. Consider the isobaric representation
$\pi=\tau \boxplus \tau_2 \boxplus \cdots \boxplus \tau_r$
of $\GL_{2m+1}(\BA)$, where $2m+1=a+\sum_{i=2}^r a_i = 2n+1-2ab$. It follows that $\pi$ has central character
$\omega_{\pi}=\omega_{\tau} \cdot \prod_{i=2}^r \omega_{\tau_i} = 1$ and
$
a\leq 2m+1= 2n+1-2ab.
$

By \cite[Theorem 3.1]{GRS11}, $\pi$ descends to an irreducible generic cuspidal representation $\sigma$ of
$\Sp_{2n-2ab}(\BA)$, which has the functorial transfer back to $\pi$. As remarked before, this is Part (3) of
Conjecture \ref{cubmfc} for the generic global Arthur parameter
$$
\psi_\pi=(\tau,1) \boxplus (\tau_2,1) \boxplus \cdots \boxplus (\tau_r,1).
$$
Hence $L(s, \tau \times \sigma)$ has a (simple) pole at $s=1$.

Let $\Delta(\tau, b)$ be the Speh residual representation in the discrete spectrum of
$\GL_{ab}(\BA)$ (see \cite{MW89}, or \cite[Section 1.2]{JLZ13}).
For any automorphic form
$$
\phi \in \CA(N_{ab}(\BA)M_{ab}(F) \bs \Sp_{2ab+2m}(\BA))_{\Delta(\tau,b) \otimes \sigma},
$$
following \cite{L76} and \cite{MW95}, one has a residual Eisenstein series
$$
E(\phi,s)(g)=E(g,\phi_{\Delta(\tau,b) \otimes \sigma},s).
$$
We refer \cite{JLZ13} for particular details about this family of Eisenstein series. In particular, it is proved in \cite{JLZ13} that
$E(\phi,s)(g)$ has a simple pole at $\frac{b+1}{2}$, which is the right-most one. We denote by $\CE(g,\phi)$ the residue,
which is square-integrable. They generate the residual representation $\CE_{\Delta(\tau, b)\otimes \sigma}$ of $\Sp_{2n}(\BA)$.
Following \cite[Section 6.2]{JLZ13}, the global Arthur parameter of this nonzero square-integrable automorphic representation
$\CE_{\Delta(\tau,b)\otimes \sigma}$ is exactly
$
\psi=(\tau, 2b+1) \boxplus \boxplus_{i=2}^r (\tau_i,1)
$
as in \eqref{case1psi}. We prove Part (3) of Conjecture \ref{cubmfc} for
\textbf{Case \Rmnum{1}}.

\begin{thm}\label{main1}
For any global Arthur parameter of the form
$$
\psi=(\tau, 2b+1) \boxplus \boxplus_{i=2}^r (\tau_i,1)
$$
with $b \geq 1$, and $\tau \ncong \tau_i$ for any $2 \leq i \leq r$,
the residual representation $\CE_{\Delta(\tau, b)\otimes \sigma}$ has a nonzero Fourier coefficient attached to
the Barbasch-Vogan duality $\eta_{\frak{so}_{2n+1}, \frak{sp}_{2n}}(\ul{p}(\psi))$ of the partition $\ul{p}(\psi)$ associated to
$(\psi,\SO_{2n+1}(\BC))$.
\end{thm}

In order to prove Theorem \ref{main1}, we have to figure out the partition $\eta_{\frak{so}_{2n+1}, \frak{sp}_{2n}}(\ul{p}(\psi))$ precisely.
We recall

\begin{defn}\label{def1}
Given any partition $\ul{q} =[q_1 q_2 \cdots q_r]$ for $\frak{so}_{2n+1}(\BC)$ with $q_1 \geq q_2 \geq \cdots \geq q_r > 0$,
whose even parts occurring with even multiplicity. Let $\ul{q}^- =[q_1 q_2 \cdots q_{r-1} (q_r-1)]$. Then the Barbasch-Vogan duality
$\eta_{\frak{so}_{2n+1}, \frak{sp}_{2n}}$, following \cite[Definition A1]{BV85} and \cite[Section 3.5]{Ac03}, is defined by
$$
\eta_{\frak{so}_{2n+1}, \frak{sp}_{2n}}
(\ul{q}) := ((\ul{q}^-)_{\Sp_{2n}})^t,
$$
where $(\ul{q}^-)_{\Sp_{2n}}$ is the $\Sp_{2n}$-collapse of $\ul{q}^-$, which is the biggest special symplectic partition which is smaller than $\ul{q}^-$.
\end{defn}

Following \cite[Section 4]{J14}, one has
$
\ul{p}(\psi)=[(2b+1)^{a}(1)^{2m+1-a}].
$
As calculated in \cite{JL15b}, when $a=2m+1$, by Definition \ref{def1}, one has
$$
\eta_{\frak{so}_{2n+1}, \frak{sp}_{2n}}(\ul{p}(\psi))
= [(a)^{2b}(2m)];
$$
when $a \leq 2m$ and $a$ is even, one obtains that
$$
\eta_{\frak{so}_{2n+1}, \frak{sp}_{2n}}(\ul{p}(\psi))
=  [(2m)(a)^{2b}];
$$
and finally, when $a \leq 2m$ and $a$ is odd, one obtains that
$$
\eta_{\frak{so}_{2n+1}, \frak{sp}_{2n}}(\ul{p}(\psi))
= [(2m)(a+1)(a)^{2b-2}(a-1)].
$$

The proof of Theorem \ref{main1} goes as follows.
Given a symplectic partition $\ul{p}$ of $2n$ (that is, odd parts occur with even multiplicities),
Denote by $\ul{p}^{\Sp_{2n}}$ the $\Sp_{2n}$-expansion of $\ul{p}$, which is the smallest special symplectic partition that is
bigger than $\ul{p}$.
In \cite{JL15a}, we proved the following theorem which provides a crucial reduction in the proof of Theorem \ref{main1}.

\begin{thm}[Theorem 4.1 \cite{JL15a}]\label{special}
Let $\pi$ be an irreducible automorphic representation of
$\Sp_{2n}(\BA)$. If $\pi$ has a nonzero Fourier coefficient attached to a
non-special symplectic partition $\ul{p}$ of $2n$, then $\pi$ must have a nonzero
Fourier coefficient attached to $\ul{p}^{\Sp_{2n}}$, the $\Sp_{2n}$-expansion of the partition $\ul{p}$.
\end{thm}

If $a \leq 2m$ and $a$ is odd, by \cite[Lemma 6.3.9]{CM93}, one has that
$$[(2m)(a+1)(a)^{2b-2}(a-1)] = [(2m)(a)^{2b}]^{\Sp_{2n}}.$$
Hence it suffices to prove the following theorem.

\begin{thm}\label{main12}
With notation above, the following hold.
\begin{enumerate}
\item If $a=2m+1$, then $\CE_{\Delta(\tau, b)\otimes \sigma}$ has a nonzero Fourier coefficient attached to $[(a)^{2b}(2m)]$.
\item If $a \leq 2m$, then $\CE_{\Delta(\tau, b)\otimes \sigma}$ has a nonzero Fourier coefficient attached to $[(2m)(a)^{2b}]$.
\end{enumerate}
\end{thm}

Part (1) and Part (2) of Theorem \ref{main12} will be proved in Sections 4 and 5, respectively.

\subsection{\textbf{Case \Rmnum{2}}}
$\psi\in\wt{\Psi}_2(\Sp_{2n})$ is written as
\begin{equation}\label{case2psi}
\psi=(\tau, 2b+1) \boxplus (\tau, 1) \boxplus \boxplus_{i=3}^r (\tau_i,1),
\end{equation}
where $b \geq 1$, $\tau \ncong \tau_i$, for any $3 \leq i \leq r$. Assume $\tau\in\CA_\cusp(\GL_a)$ has central character $\omega_{\tau}$,
and $\tau_i\in\CA_\cusp(\GL_{a_i})$ has central character $\omega_{\tau_i}$, $3 \leq i \leq r$.
Then $2n+1=a(2b+1)+ a + \sum_{i=3}^r {a_i}$ and $\omega_{\tau}^{2b+1} \cdot \omega_{\tau} \cdot \prod_{i=3}^r \omega_{\tau_i} = 1$.
Consider the isobaric representation
$\pi=\tau_3 \boxplus \cdots \boxplus \tau_r$
of $\GL_{2m+1}(\BA)$, where $2m+1=\sum_{i=3}^r a_i = 2n+1-a(2b+2)$. Then $\pi$ has central character
$\omega_{\pi}=\prod_{i=3}^r \omega_{\tau_i} = 1$.

By \cite[Theorem 3.1]{GRS11}, there is a generic $\sigma\in\CA_\cusp(\Sp_{2m})$
such that $\sigma$ has the functorial transfer $\pi$ and hence $L(s, \tau \times \sigma)$ is holomorphic at $s=1$ in this case.
For any automorphic form
$$
\phi \in \CA(N_{a(b+1)}(\BA)M_{a(b+1)}(F) \bs \Sp_{2a(b+1)+2m}(\BA))_{\Delta(\tau,b+1) \otimes \sigma},
$$
one defines a residual Eisenstein series as in \textbf{Case \Rmnum{1}}
$$
E(\phi,s)(g)=E(g,\phi_{\Delta(\tau,b+1) \otimes \sigma},s)
$$
By \cite{JLZ13}, this Eisenstein series has a simple pole at $\frac{b}{2}$, which is the right-most one.
Denote the representation generated by these residues at $s=\frac{b}{2}$
by $\CE_{\Delta(\tau, b+1)\otimes \sigma}$, which is square-integrable.
Following \cite{JLZ13} and \cite[Theorem 7.1.2]{Sh10}, this residual representation $\CE_{\Delta(\tau, b+1)\otimes \sigma}$ is nonzero.
In particular, by Section 6.2 of \cite{JLZ13}, the global Arthur parameter of $\CE_{\Delta(\tau,b+1)\otimes \sigma}$ is exactly
$
\psi=(\tau, 2b+1) \boxplus (\tau,1) \boxplus \boxplus_{i=3}^r (\tau_i,1)
$
as in \textbf{Case \Rmnum{2}}.
In this case, we prove

\begin{thm}\label{main2}
For any global Arthur parameter of the form
$$
\psi=(\tau, 2b+1) \boxplus (\tau, 1) \boxplus \boxplus_{i=3}^r (\tau_i,1)
$$
with $b \geq 1$, and $\tau \ncong \tau_i$ for any $3 \leq i \leq r$,
the residual representation $\CE_{\Delta(\tau, b+1)\otimes \sigma}$
has a nonzero Fourier coefficient attached to
the Barbasch-Vogan duality $\eta_{\frak{so}_{2n+1}, \frak{sp}_{2n}}(\ul{p}(\psi))$ of the partition $\ul{p}(\psi)$ associated to
$(\psi,\SO_{2n+1}(\BC))$.
\end{thm}

Following \cite[Section 4]{J14}, one has
$
\ul{p}(\psi)=[(2b+1)^{a}(1)^a(1)^{2m+1}].
$
By Definition \ref{def1}, we may calculate the partition $\eta_{\frak{so}_{2n+1}, \frak{sp}_{2n}}(\ul{p}(\psi))$ explicitly as follows.
When $a$ is even, we have
\begin{align*}
\eta_{\frak{so}_{2n+1}, \frak{sp}_{2n}}(\ul{p}(\psi))
= \ & \eta_{\frak{so}_{2n+1}, \frak{sp}_{2n}}([(2b+1)^{a}(1)^{2m+1+a}])\\
= \ & [(2b+1)^a(1)^{2m+a}]^t\\
= \ & [(a)^{2b+1}] + [(2m+a)]\\
= \ & [(2m+2a)(a)^{2b}].
\end{align*}
When $a$ is odd, we have
\begin{align*}
\eta_{\frak{so}_{2n+1}, \frak{sp}_{2n}}(\ul{p}(\psi))
= \ & \eta_{\frak{so}_{2n+1}, \frak{sp}_{2n}}([(2b+1)^{a}(1)^{2m+1+a}])\\
= \ & ([(2b+1)^a(1)^{2m+a}]_{\Sp_{2n}})^t\\
= \ & [(2b+1)^{a-1}(2b)(2)(1)^{2m+a-1}]^t\\
= \ & [(a-1)^{2b+1}] + [(1)^{2b}] + [(1)^2] + [(2m+a-1)]\\
= \ & [(2m+2a)(a+1)(a)^{2b-2}(a-1)].
\end{align*}
As before, if $a$ is odd, then, by the recipe of
obtaining $\Sp_{2n}$-expansion of a symplectic partition $\ul{p}$
given in \cite[Lemma 6.3.9]{CM93}, we have that
$$
[(2m+2a)(a+1)(a)^{2b-2}(a-1)] = [(2m+2a)(a)^{2b}]^{\Sp_{2n}}.
$$
Hence it suffices to prove the following theorem.

\begin{thm}\label{main22}
The residual representation $\CE_{\Delta(\tau, b+1)\otimes \sigma}$
has a nonzero Fourier coefficient attached to $[(2m+2a)(a)^{2b}]$.
\end{thm}

The proof of Theorem \ref{main22} is given in Section 6, using induction on the integer $b$. We note that when $b=0$, the
Arthur parameter is
$$
\psi=2(\tau, 1) \boxplus \boxplus_{i=3}^r (\tau_i,1),
$$
which does not parameterize automorphic representations in the discrete spectrum. Indeed, in this case, the corresponding automorphic
representation constructed from the Eisenstein series is the value at $s=0$, which we still denote it by
$\CE_{\Delta(\tau,1)\otimes\sigma}=\CE_{\tau\otimes\sigma}$. It is clear that in this case, the partition $\ul{p}(\psi)$ is the trivial
partition. On the other hand, following \cite[Theorem 7.1.3]{Sh10}, the representation $\CE_{\Delta(\tau,1)\otimes\sigma}$ has a
nonzero Whittaker-Fourier coefficient. In other words, Theorem \ref{main22} still holds for $b=0$. As we proceed in Section 6, the
case of $b=0$ will serve as the base of the induction argument.

\subsection{\textbf{Case \Rmnum{3}}}
$\psi\in\wt{\Psi}_2(\Sp_{2n})$ is written as
\begin{equation}\label{case3psi}
\psi=(\tau, 2b) \boxplus \boxplus_{i=2}^r (\tau_i,1),
\end{equation}
where $b \geq 1$. In this case, $\tau$ is of symplectic type (and hence $a=2k$ is even),
while $\tau_i$ for all $2 \leq i \leq r$ are of orthogonal type.
Assume that $\tau\in\CA_\cusp(\GL_a)$ has central character $\omega_{\tau}$, and
$\tau_i\in\CA_\cusp(\GL_{a_i})$ has central character $\omega_{\tau_i}$, $2 \leq i \leq r$.
By the definition of Arthur parameters, one has that $2n+1=2ab + \sum_{i=2}^r {a_i}$, and $\prod_{i=2}^r \omega_{\tau_i} = 1$.
Consider the isobaric representation
$\pi=\tau_2 \boxplus \cdots \boxplus \tau_r$
of $\GL_{2m+1}(\BA)$, where $2m+1=\sum_{i=2}^r a_i = 2n+1-2ab$. Hence $\pi$ has central character
$\omega_{\pi}=\prod_{i=2}^r \omega_{\tau_i} = 1$.

By \cite[Theorem 3.1]{GRS11}, there is a generic $\sigma\in\CA_\cusp(\Sp_{2m})$ that
has the functorial transfer $\pi$. Then we define a residual Eisenstein series
$$
E(,\phi,s)(g)=E(g,\phi_{\Delta(\tau,b) \otimes \sigma},s)
$$
associated to any automorphic form
$$
\phi \in \CA(N_{ab}(\BA)M_{ab}(F) \bs \Sp_{2ab+2m}(\BA))_{\Delta(\tau,b) \otimes \sigma}.
$$
By \cite{JLZ13}, this Eisenstein series
may have a simple pole at $\frac{b}{2}$, which is the right-most one.
Denote the representation generated by these residues at $s=\frac{b}{2}$
by $\CE_{\Delta(\tau, b)\otimes \sigma}$.
This residual representation is square-integrable. If $L(\frac{1}{2}, \tau \times \sigma) \neq 0$, the residual
representation $\CE_{\tau\otimes \sigma}$ is nonzero, and hence by the induction argument in \cite{JLZ13}, the
residual representation $\CE_{\Delta(\tau, b)\otimes \sigma}$ is also nonzero.
Finally, following Section 6.2 of \cite{JLZ13}, we see that the global Arthur parameter of $\CE_{\Delta(\tau,b)\otimes \sigma}$ is exactly
$
\psi=(\tau, 2b)\boxplus \boxplus_{i=2}^r (\tau_i,1)
$
as in \eqref{case3psi}. We prove

\begin{thm}\label{main3}
Assume that $a=2k$ and $L(\frac{1}{2}, \tau \times \sigma) \neq 0$.
If the residual representation $\CE_{\tau\otimes \sigma}$ of $\Sp_{4k+2m}(\BA)$, with $\sigma \ncong 1_{\Sp_0(\BA)}$,
has a nonzero Fourier coefficient attached to the partition $[(2k+2m)(2k)]$, then, for any $b \geq 1$, the residual
representation $\CE_{\Delta(\tau, b)\otimes \sigma}$ has a nonzero Fourier coefficient attached to the partition
$[(2k+2m)(2k)^{2b-1}]$.
\end{thm}

We remark that if $\sigma \cong 1_{\Sp_0(\BA)}$,
$L(\frac{1}{2}, \tau)=L(\frac{1}{2}, \tau \times \sigma) \neq 0$. In this case, \cite[Theorem 4.2.2]{L13b} shows that
$\mathfrak{p}^m(\CE_{\Delta(\tau, b)\otimes \sigma}) = \{[(2k)^{2b}]\}$.

In fact, the assumption that the residual representation $\CE_{\tau\otimes \sigma}$ of $\Sp_{4k+2m}(\BA)$,
with $\sigma \ncong 1_{\Sp_0(\BA)}$, has a nonzero Fourier coefficient attached to the partition $[(2k+2m)(2k)]$
is exactly \cite[Conjecture 6.1]{GJR04}, and hence Theorem \ref{main3} has a close connection to the Gan-Gross-Prasad conjecture (\cite{GGP12}).
We will come back to this issue in our future work.

In this case, one has that $\ul{p}(\psi)=[(2b)^{a}(1)^{2m+1}]$, and following the calculation in \cite{JL15b}, one has
$$
\eta_{\frak{so}_{2n+1}, \frak{sp}_{2n}}(\ul{p}(\psi))
= [(a+2m)(a)^{2b-1}],
$$
where $a=2k$ is even. The proof of Theorem \ref{main3} is given in Section 7.

When $L(\frac{1}{2}, \tau \times \sigma)$ is zero for
the Arthur parameter in \eqref{case3psi}, the corresponding automorphic $L^2$-packet $\wt{\Pi}_\psi(\epsilon_\psi)$
are expected to contain all cuspidal automorphic representations if it is not empty. We are going to apply the construction of
endoscopy correspondences outlined in \cite{J14} to construct the distinguished cuspidal members in $\wt{\Pi}_\psi(\epsilon_\psi)$.
The details for this case will be considered in our future work. See \cite{JL15b} for a brief discussion in this aspect.


\section{A basic lemma}


We recall a basic lemma from \cite{JL15b}, which will be a technical key step in the proofs of this paper.
Let $H$ be a reductive group defined over $F$.
We first recall \cite[Lemma 5.2]{JL13}, which is also formulated in a slightly different version in \cite[Corollary 7.1]{GRS11}. Note that the proof of \cite[Lemma 5.2]{JL13} is valid for $H(\BA)$.

Let $C$ be an $F$-subgroup of a maximal unipotent subgroup of $H$, and let $\psi_C$ be a non-trivial character of $[C] = C(F) \bs C(\BA)$.
$\wt{X}, \wt{Y}$ are two unipotent $F$-subgroups, satisfying the following conditions:
\begin{itemize}
\item[(1)] $\wt{X}$ and $\wt{Y}$ normalize $C$;
\item[(2)] $\wt{X} \cap C$ and $\wt{Y} \cap C$ are normal in $\wt{X}$ and $\wt{Y}$, respectively, $(\wt{X} \cap C) \bs \wt{X}$ and $(\wt{Y} \cap C) \bs \wt{Y}$ are abelian;
\item[(3)] $\wt{X}(\BA)$ and $\wt{Y}(\BA)$ preserve $\psi_C$;
\item[(4)] $\psi_C$ is trivial on $(\wt{X} \cap C)(\BA)$ and $(\wt{Y} \cap C)(\BA)$;
\item[(5)] $[\wt{X}, \wt{Y}] \subset C$;
\item[(6)]  there is a non-degenerate pairing $(\wt{X} \cap C)(\BA) \times (\wt{Y} \cap C)(\BA) \rightarrow \BC^*$, given by $(x,y) \mapsto \psi_C([x,y])$, which is
multiplicative in each coordinate, and identifies $(\wt{Y} \cap C)(F) \bs \wt{Y}(F)$ with the dual of
$
\wt{X}(F)(\wt{X} \cap C)(\BA) \bs \wt{X}(\BA),
$
and
$(\wt{X} \cap C)(F) \bs \wt{X}(F)$ with the dual of
$
\wt{Y}(F)(\wt{Y} \cap C)(\BA) \bs \wt{Y}(\BA).
$
\end{itemize}

Let $B =C\wt{Y}$ and $D=C\wt{X}$, and extend $\psi_C$ trivially to characters of $[B]=B(F)\bs B(\BA)$ and $[D]=D(F)\bs D(\BA)$,
which will be denoted by $\psi_B$ and $\psi_D$ respectively.

\begin{lem}[Lemma 5.2 of \cite{JL13}]\label{nvequ}
Assume that $(C, \psi_C, \wt{X}, \wt{Y})$ satisfies all the above conditions. Let $f$ be an automorphic form on $H(\BA)$. Then
$$\int_{[C]} f(cg) \psi_C^{-1}(c) dc \equiv 0, \forall g \in H(\BA),$$
if and only if
$$\int_{[D]} f(ug) \psi_D^{-1}(u) du \equiv 0, \forall g \in H(\BA),$$
if and only if
$$\int_{[B]} f(vg) \psi_B^{-1}(v) dv \equiv 0, \forall g \in H(\BA).$$
\end{lem}

For simplicity, we always use $\psi_C$ to denote its extensions $\psi_B$ and $\psi_D$ when we apply Lemma \ref{nvequ} to various circumstances.
Lemma \ref{nvequ} can be extended as follows and will be a technical key in this paper.

\begin{lem}[Lemma 6.2 of \cite{JL15b}]\label{nvequ2}
Assume the quadruple $(C, \psi_C, \wt{X}, \wt{Y})$ satisfies the following conditions:
$\wt{X} = \{\wt{X}_i\}_{i=1}^r$, $\wt{Y} = \{\wt{Y}_i\}_{i=1}^r$, for $1 \leq i \leq r$, each quadruple
$$(\wt{X}_{i-1} \cdots \wt{X}_1 C \wt{Y}_r \cdots \wt{Y}_{i+1}, \psi_C, \wt{X}_i, \wt{Y}_i)$$
satisfies all the conditions for Lemma \ref{nvequ}.
Let $f$ be an automorphic form on $H(\BA)$.
Then
$$\int_{[\wt{X}_{r} \cdots \wt{X}_1 C]} f(xcg) \psi_C^{-1}(c) dcdx \equiv 0, \forall g \in H(\BA),$$
if and only if
$$\int_{[C\wt{Y}_r \cdots \wt{Y}_{1}]} f(cyg) \psi_C^{-1}(c) dydc \equiv 0, \forall g \in H(\BA).$$
\end{lem}

The proof of this lemma is carried out by using Lemma \ref{nvequ} inductively, and was given with full details in \cite{JL15b}.

\section{Proof of Part (1) of Theorem \ref{main12}}


In this section, we assume that $a=2m+1$ and show that $\CE_{\Delta(\tau, b)\otimes \sigma}$ has a nonzero Fourier coefficient attached to $\ul{p} := [(2m+1)^{2b}(2m)]$.

\textbf{Proof of Part (1) of Theorem \ref{main12}}.
We will prove by induction on $b$. Note that when $b=0$,
$\CE_{\Delta(\tau, b)\otimes \sigma} \cong \sigma$ which has a nonzero Fourier coefficient attached to $[(2m)]$ since $\sigma$ is generic.
Assume that $\CE_{\Delta(\tau, b-1)\otimes \sigma}$ has a nonzero $\psi_{[(2m+1)^{2b-2}(2m)], \alpha}$-Fourier coefficient attached to $[(2m+1)^{2b-2}(2m)]$,
for some $\alpha \in F^*/ (F^*)^2$.

Take any $\varphi \in \CE_{\Delta(\tau, b)\otimes \sigma}$, consider its $\psi_{\ul{p}, \alpha}$-Fourier coefficients
attached to $\ul{p}$ as follows
\begin{equation}\label{part1equ1}
\varphi^{\psi_{\ul{p}, \alpha}}(g) = \int_{[V_{\ul{p},2}]}
\varphi(vg) \psi_{\ul{p}, \alpha}^{-1}(v)dv.
\end{equation}
For definitions of the unipotent group $V_{\ul{p},2}$
and its character $\psi_{\ul{p}, \alpha}$, see \cite[Section 2]{JL15a}. By \cite[Corollary 2.4]{JL15a}, the integral in
\eqref{part1equ1} is non-vanishing if and only if the following integral is non-vanishing:
\begin{equation}\label{part1equ2}
\int_{[Y_1 V_{\ul{p},2}]}
\varphi(vg) \psi_{\ul{p}, \alpha}^{-1}(v)dv,
\end{equation}
where $Y_1$ is defined in \cite[(2.5)]{JL15a} corresponding to the partition $[(2m+1)^{2b}(2m)]$, and the character $\psi_{\ul{p}, \alpha}$ extends to $Y_1 V_{\ul{p},2}$ trivially.

Assume that $T$ is the maximal split torus in $\Sp_{2b(2m+1)+2m}$,
consisting of elements
$$
\diag(t_1, t_2, \ldots, t_{b(2m+1)+m}, t_{b(2m+1)+m}^{-1}, \ldots, t_2^{-1}, t_1^{-1}).
$$
Let $\omega_1$ be the Weyl element of $\Sp_{2b(2m+1)+2m}$, sending elements $t \in T$
to the following torus elements:
\begin{equation}\label{part1equ3}
t'=\diag(t^{(0)}, t^{(1)}, t^{(2)}, \ldots, t^{(m)}, t^{(m+1)}, t^{(m),*}, \ldots, t^{(2),*}, t^{(1),*}, t^{(0),*}),
\end{equation}
where $t^{(0)} =  \diag(t_1, t_2, \ldots, t_{2m+1})$, and with $e=2m+1$,
$$
t^{(m+1)} =  \diag(t_{e+m+1}, \ldots, t_{(b-1)e+m+1}, t_{be-m}^{-1}, \ldots, t_{2e-m}^{-1})
$$
and
$$
t^{(j)} =  \diag(t_{e+j}, \ldots, t_{(b-1)e+j}, t_{be-j+1}^{-1}, \ldots, t_{2e-j+1}^{-1}, t_{be+j}),
$$
for $1 \leq j \leq m$.

Identify $\Sp_{(2b-1)(2m+1)+2m}$ with its image in $\Sp_{2b(2m+1)+2m}$ under
the embedding $g \mapsto \diag(I_{2m+1}, g, I_{2m+1})$. Denote the restriction of $\omega_1$ to $\Sp_{(2b-1)(2m+1)+2m}$
by $\omega_1'$. We conjugate cross the integration variables by $\omega_1$ from the left and obtain that the integral in \ref{part1equ2} becomes:
\begin{equation}\label{part1equ4}
\int_{[U_{\ul{p},2}]} \varphi(u \omega_1 g) \psi_{\underline{p},\alpha}^{\omega_1}(u)^{-1} du,
\end{equation}
where $U_{\ul{p},2} = \omega_1 Y_1V_{\ul{p},2} \omega_1^{-1}$,
and $\psi_{\underline{p},\alpha}^{\omega_1}(u) = \psi_{\underline{p},\alpha}(\omega_1^{-1}u\omega_1)$.

Now, we describe the structure of elements in $U_{\ul{p},2}$.
Any element in $U_{\underline{p},2}$
has the following form:
\begin{equation}\label{part1equ5}
u =
\begin{pmatrix}
z_{2m+1} & q_1 & q_2\\
0 & u' & q_1^*\\
0 & 0 & z_{2m+1}^*
\end{pmatrix}
\begin{pmatrix}
I_{2m+1} & 0 & 0\\
p_1 & I_{((2b-2)(2m+1)+2m} & 0\\
p_2 & p_1^* & I_{2m+1}
\end{pmatrix},
\end{equation}
where $z_{2m+1} \in V_{2m+1}$, the standard maximal unipotent subgroup of $\GL_{2m+1}$;
$u' \in U_{[(2m+1)^{2b-2}(2m)],2} := \omega_1' Y_2V_{[(2m+1)^{2b-2}(2m)],2} \omega_1'^{-1}$ with $Y_2$ as in \cite[(2.5)]{JL15a}
corresponding to the partition $[(2m+1)^{2b-2}(2m)]$; and
$p_i, q_i$, $1 \leq i \leq 2$, are described as follows.
\begin{itemize}
\item $q_1 \in M_{(2m+1)\times ((2b-2)(2m+1)+2m)}$, such that $q_1(i,j)=0$, for $1 \leq i \leq 2m+1$
and $1 \leq j \leq (2b-2)+(2b-1)(i-1)$.
\item $p_1 \in M_{((2b-2)(2m+1)+2m) \times (2m+1)}$,
such that $p_1(i,j) = 0$, for $1 \leq j \leq 2m+1$
and $(2b-2)+(2b-1)(i-1) + 1 \leq i \le (2b-2)(2m+1)+2m$.
\item $q_2 \in M_{(2m+1)\times (2m+1)}$, symmetric with respect to the secondary diagonal, such that $q_2(i,j)=0$, for
$1 \leq i \leq 2m+1$ and $1 \leq j \leq i$.
\item $p_2 \in M_{(2m+1)\times (2m+1)}$, symmetric with respect to the secondary diagonal, such that $p_2(i,j)=0$, for
$1 \leq i \leq 2m+1$ and $1 \leq j \leq i$.
\end{itemize}
Note that
\begin{align*}
& \psi_{\underline{p},\alpha}^{\omega_1}
(\begin{pmatrix}
z_{2m+1} & q_1 & q_2\\
0 & I_{(2b-2)(2m+1)+2m} & q_1^*\\
0 & 0 & z_{2m+1}^*
\end{pmatrix})
=  \psi(\sum_{i=1}^{2m} z_{2m+1}(i,i+1)).
\end{align*}

Next, we apply Lemma \ref{nvequ2} to fill the zero entries in $q_1, q_2$ using the nonzero entries in $p_1, p_2$.
To proceed, we need to define a sequence of one-dimensional root subgroups and put them in a correct order.

Let $X_j$'s, with $1 \leq j \leq (2b-2)+1$, be the one-dimensional subgroups
corresponding to the roots such that
the corresponding entries are in the first row of $q_1$ or $q_2$ and are identically zero, from right to left.
For $1 < i \leq m$, let $X_j$'s, with
$$
(\sum_{k=1}^{i-1} [(2b-2)+(2b-1)(k-1)+k]) +1 \leq j \leq \sum_{k=1}^{i} [(2b-2)+(2b-1)(k-1)+k],
$$
be the one-dimensional subgroups corresponding to the roots such that the corresponding entries are in the $i$-th row of $q_1$ or $q_2$ and are identically zero, from right to left.

Let $Y_j$'s, with $1 \leq j \leq (2b-2)+1$, be the one-dimensional subgroups corresponding to the roots such that the corresponding entries are in the second column of $p_1$ or $p_2$ and are not identically zero, from bottom to top.
For $1 < i \leq m$, let $Y_j$'s, with
$$
1+\sum_{k=1}^{i-1} [(2b-2)+(2b-1)(k-1)+k] \leq j \leq \sum_{k=1}^{i} [(2b-2)+(2b-1)(k-1)+k],
$$
be the one-dimensional subgroups corresponding to the roots such that the corresponding entries are in the $(i+1)$-th column of $p_1$ or $p_2$ and are not identically zero, from bottom to top.

Let $W_1$ be the subgroup of $U_{\ul{p}, 2}$ such that the entries corresponding to the one-dimensional subgroups $Y_j$'s above, with
$$
1 \leq j \leq \ell := \sum_{k=1}^{m} [(2b-2)+(2b-1)(k-1)+k],
$$
are all identically zero. And let $\psi_{W_1} = \psi_{\underline{p},\alpha}^{\omega_1}|_{W_1}$.
Then the quadruple
$(W_1, \psi_{W_1}, \{X_j\}_j^{\ell}, \{Y_j\}_j^{\ell})$ satisfies all the conditions for Lemma \ref{nvequ2}. Hence, by Lemma \ref{nvequ2},
the integral in \eqref{part1equ4} is non-vanishing if and only if the following integral is non-vanishing:
\begin{equation}\label{part1equ6}
\int_{[W_2]} \varphi(w \omega_1 g) \psi_{W_2}(w)^{-1} dw,
\end{equation}
where $W_2 := \prod_{j=1}^{\ell} X_j W_1$, and
$\psi_{W_2}$ is the character on $W_2$ extended trivally from
$\psi_{W_1}$.

Now we consider the $i$-th row of $q_1$ and $q_2$, $m+1 \leq i \leq 2m$. We will continue to apply Lemma \ref{nvequ2} to fill the zero entries in $q_1$ and $q_2$, row-by-row, from the $m+1$-th row to $2m$-th row. But for each $m+1 \leq i \leq 2m$, before we apply Lemma \ref{nvequ2} as above, we need to take the Fourier expansion along the one-dimensional root subgroup $X_{2e_i}$. For example, for $i=m+1$.
we first take the Fourier expansion of the integral in \eqref{part1equ6} along the one-dimensional root subgroup $X_{2e_{m+1}}$. We will get two kinds of Fourier coefficients corresponding to the orbits of the dual of $[X_{2e_{m+1}}]:=X_{2e_{m+1}}(F) \bs X_{2e_{m+1}}(\BA)$: the trivial orbit and the non-trivial one. For the Fourier coefficients attached to the non-trivial orbit,
we can see that there is an inner integral $\varphi^{\psi_{[(2m+2)1^{2b(2m+1)-2}],\beta}}$
($\beta \in F^*$), which is identically zero by \cite[Proposition 6.4]{JL15c}.
Therefore only the Fourier coefficient attached to the trivial orbit,
which actually equals to the integral in \eqref{part1equ6},
survives. Then, we can apply the Lemma \ref{nvequ2} to the $m+1$-th row of $q_1$ and $q_2$ similarly as above.

After considering all the $i$-th row of $q_1$ and $q_2$, $m+1 \leq i \leq 2m$ as above, we get that the integral in \eqref{part1equ6} is non-vanishing if and only if the following integral is non-vanishing:
\begin{equation}\label{part1equ7}
\int_{[W_3]} \varphi(w \omega_1 g) \psi_{W_3}(w)^{-1} dw,
\end{equation}
where $W_3$ has elements of the following form:

\begin{equation}\label{part1equ8}
w =
\begin{pmatrix}
z_{2m+1} & q_1 & q_2\\
0 & u' & q_1^*\\
0 & 0 & z_{2m+1}^*
\end{pmatrix},
\end{equation}
where $z_{2m+1} \in V_{2m+1}$, the standard maximal unipotent subgroup of $\GL_{2m+1}$;
$u' \in U_{[(2m+1)^{2b-2}(2m)],2} := \omega_1' Y_2V_{[(2m+1)^{2b-2}(2m)],2} \omega_1'^{-1}$ with $Y_2$ as in \cite[(2.5)]{JL15a}
corresponding to the partition $[(2m+1)^{2b-2}(2m)]$;
$q_1 \in M_{(2m+1)\times ((2b-2)(2m+1)+2m)}$, such that $q_1(2m+1,j)=0$, for $1 \leq j \leq (2b-2)(2m+1)+2m$;
$q_2 \in M_{(2m+1)\times (2m+1)}$, symmetric with respect to the secondary diagonal, such that $q_2(2m+1,1)=0$.
And
\begin{align*}
& \psi_{W_3}
(\begin{pmatrix}
z_{2m+1}  & q_1 & q_2\\
0 & I_{(2b-2)(2m+1)+2m} & q_1^*\\
0 & 0 & z_{2m+1} ^*
\end{pmatrix})
=  \psi(\sum_{i=1}^{2m} z_{2m+1}(i,i+1)).
\end{align*}

Now consider the Fourier expansion of the integral in \eqref{part1equ7} along the one-dimensional root subgroup $X_{2e_{2m+1}}$. By the same reason as above, only the Fourier coefficient corresponding to the trivial orbit of the dual of
$[X_{2e_{2m+1}}]$ survives, which is actually equal to the integral in \eqref{part1equ7}:
\begin{equation}\label{part1equ9}
\int_{[W_4]} \varphi(w \omega_1 g) \psi_{W_4}(w)^{-1} dw,
\end{equation}
where elements in $W_4$ have the same structure as in \eqref{part1equ8}, except that $q_2(2m+1,1)$ is not identically zero.

It is easy to see that the integral in \eqref{part1equ9} has an inner integral which is exactly $\varphi^{\psi_{N_{1^{2m}}}}$, using notation in Lemma \ref{lem2} below. On the other hand, we know that by Lemma \ref{lem2} below,
$\varphi^{\psi_{N_{1^{2m}}}} = \varphi^{\wt{\psi}_{N_{1^{2m+1}}}}$.
Therefore, the
integral in \eqref{part1equ9}
becomes
\begin{equation}\label{part1equ10}
\int_{[W_5]} \varphi(w \omega_1 g) \psi_{W_5}(w)^{-1} dw,
\end{equation}
where elements in $W_5$ are of the form:
$$
w = w(z_{2m+1}, u',q_1, q_2) =
\begin{pmatrix}
z_{2m+1} & q_1 & q_2\\
0 & u' & q_1^*\\
0 & 0 & z_{2m+1}^*
\end{pmatrix},
$$
where $z_{2m+1} \in V_{2m+1}$, the standard maximal unipotent subgroup of
$\GL_{2m+1}$, $u' \in U_{[(2m+1)^{2b-2}(2m)],2} := \omega_1' Y_2V_{[(2m+1)^{2b-2}(2m)],2} \omega_1'^{-1}$ with $Y_2$ as in
\cite[(2.5)]{JL15a} corresponding to the partition $[(2m+1)^{2b-2}(2m)]$; $q_1 \in M_{(2m+1)\times ((2b-2)(2m+1)+2m)}$, and
$q_2 \in M_{(2m+1)\times (2m+1)}$, symmetric with respect to the secondary diagonal.
And
\begin{align*}
& \psi_{W_5}
(\begin{pmatrix}
z_{2m+1} & q_1 & q_2\\
0 & I_{(2b-2)(2m+1)+2m} & q_1^*\\
0 & 0 & z_{2m+1}^*
\end{pmatrix})
=  \psi(\sum_{i=1}^{2m} z_{2m+1}(i,i+1)).
\end{align*}
Hence, the integral in \eqref{part1equ10} can be written as
\begin{equation}\label{part1equ11}
\int_{W_6} \varphi_{P_{2m+1}}(w \omega_1 g) \psi_{W_6}(w)^{-1} dw,
\end{equation}
where
$W_6$ is a subgroup of $W_5$ consisting of elements of the form
$w(z_{2m+1}, u',0, 0)$, $\psi_{W_6}=\psi_{W_5}|_{W_6}$,
and
$\varphi_{P_{2m+1}}$ is the constant term of
$\varphi$ along the parabolic subgroup $P_{2m+1}=M_{2m+1}N_{2m+1}$ of $\Sp_{2b(2m+1)+2m}$ with the
Levi subgroup isomorphic to $\GL_{2m+1} \times \Sp_{(2b-2)(2m+1)+2m}$.

By Lemma \ref{constantterm} below,
$\varphi_{P_{2m+1}}(w \omega_1 g)$ is an automorphic form
in $\tau \lvert \cdot \rvert^{-b}
\otimes \mathcal{E}_{\Delta(\tau,b-1)\otimes \sigma}$ when restricted to the Levi subgroup.
Note that the restriction of $\psi_{W_5}$
to the $z_{2m+1}$-part gives a Whittaker coefficient of $\tau$, and the restriction
to the $u'$-part gives a $\psi_{[(2m+1)^{2b-2}(2m)], \alpha}$-Fourier coefficient of $\mathcal{E}_{\Delta(\tau,b-1)\otimes \sigma}$
up to the conjugation of
the Weyl element $\omega_1'$. On the other hand, $\tau$ is
generic, and by induction assumption, $\mathcal{E}_{\Delta(\tau,b-1)\otimes \sigma}$
has a nonzero $\psi_{[(2m+1)^{2b-2}(2m)], \alpha}$-Fourier coefficient.
Therefore, we conclude that
$\mathcal{E}_{\Delta(\tau,b)\otimes \sigma}$
has a nonzero $\psi_{\underline{p},\alpha}$-Fourier coefficient
attached to the partition $\ul{p}=[(2m+1)^{2b}(2m)]$.
This completes the proof of Part (1) of Theorem \ref{main12}, up to Lemmas \ref{constantterm}
and \ref{lem2}, which are stated below.

Note that Lemmas \ref{constantterm}
and \ref{lem2} are analogues of \cite[Lemmas 4.2.4 and 4.2.6]{L13b}, with similar arguments, and hence we state them without proofs.

\begin{lem}\label{constantterm}
Let $P_{ai} = M_{ai} N_{ai}$, with $1 \leq i \leq b$ and $a \leq 2m+1$, be the parabolic subgroup of $\Sp_{2ab+2m}$ with
Levi part
$$M_{ai} \cong \GL_{ai} \times \Sp_{a(2b-2i)+2m}.
$$
Let $\varphi$
be an arbitrary automorphic form in $\CE_{\Delta(\tau, b)\otimes \sigma}$.
Denote by $\varphi_{P_{ai}}(g)$ the constant term of $\varphi$
along $P_{ai}$.
Then, for $1 \leq i \leq b$,
$$
\varphi_{P_{ai}} \in \CA(N_{ai}(\BA)M_{ai}(F)\bs \Sp_{2ab+2m}(\BA))_{\Delta(\tau, i)\lvert \cdot \rvert^{-\frac{2b+1-i}{2}}
\otimes \mathcal{E}_{\Delta(\tau,b-i) \otimes \sigma}}.
$$
Note that when $b=i$, $\mathcal{E}_{\Delta(\tau,b-i) \otimes \sigma} = \sigma$.
\end{lem}


\begin{lem}\label{lem2}
Let $N_{1^{p}}$ be the unipotent radical of the parabolic subgroup
$P_{1^{p}}$ of $\Sp_{2b(2m+1)+2m}$ with the Levi part being
$\GL_1^{\times p} \times \Sp_{2b(2m+1)+2m-2p}$. Let
$$
\psi_{N_{1^{p}}}(n) := \psi(n_{1,2}+ \cdots + n_{p,p+1}),
$$
and
$$
\widetilde{\psi}_{N_{1^{p}}}(n) := \psi(n_{1,2}+ \cdots + n_{p-1,p}),
$$
be two characters of $N_{1^{p}}$. For any automorphic form $\varphi \in \mathcal{E}_{\Delta(\tau,b)\otimes \sigma}$,
define $\psi_{N_{1^{p}}}$ and $\widetilde{\psi}_{N_{1^{p}}}$-Fourier coefficients
as follows:
\begin{equation}\label{lem2equ1}
\varphi^{\psi_{N_{1^{p}}}}(g) := \int_{[N_{1^{p}}]} \varphi(ng)\psi_{N_{1^{p}}}(n)^{-1}
dn,
\end{equation}
\begin{equation}\label{lem2equ2}
\varphi^{\widetilde{\psi}_{N_{1^{p}}}}(g) := \int_{[N_{1^{p}}]} \varphi(ng)\widetilde{\psi}_{N_{1^{p}}}(n)^{-1}
du.
\end{equation}
Then $\varphi^{\psi_{N_{1^{p}}}} \equiv 0, \forall p \geq 2m+1$, and
$\varphi^{\psi_{N_{1^{2m}}}} = \varphi^{\widetilde{\psi}_{N_{1^{2m+1}}}}$.
\end{lem}


\section{Proof of Part (2) of Theorem \ref{main12}}

In this section, we assume that $a \leq 2m$ and $\sigma$ is $\psi^{\alpha}$-generic for $\alpha \in F^*/(F^*)^2$, and show that
$\CE_{\Delta(\tau, b)\otimes \sigma}$ has a nonzero Fourier coefficient attached to $[(2m)(a)^{2b}]$.

First, we construct a residual representation of $\wt{\Sp}_{2ab}(\BA)$ as follows. For any $\wt{\phi} \in \CA(N_{ab}(\BA)\wt{M}_{ab}(F) \bs \wt{\Sp}_{2ab}(\BA))_{\gamma_{\psi^{-\alpha}} \Delta(\tau,b)}$, following \cite{MW95}, an
residual Eisenstein series
can be defined by
$$
\wt{E}(\wt{\phi},s)(g)=\sum_{\gamma\in P_{ab}(F)\bks \Sp_{2ab}(F)}\lambda_s \wt{\phi}(\gamma g).
$$
It converges absolutely for real part of $s$ large and has meromorphic continuation to the whole complex plane $\BC$. By similar argument as that in \cite{JLZ13}, this Eisenstein series
has a simple pole at $\frac{b}{2}$, which is the right-most one.
Denote the representation generated by these residues at $s=\frac{b}{2}$
by $\wt{\CE}_{\Delta(\tau, b)}$. This residual representation is square-integrable.


We separate the proof of Part (2) of Theorem \ref{main12} into following \textbf{three steps}:

\textbf{Step (1).} $\CE_{\Delta(\tau, b)\otimes \sigma}$ has a nonzero Fourier coefficient attached to the partition $[(2m)1^{2ab}]$ with respect to the character $\psi_{[(2m)1^{2ab}], \alpha}$ (for definition, see \cite[Section 2]{JL15a}).

\textbf{Step (2).} $\wt{\CE}_{\Delta(\tau,b)}$ is irreducible. Let $\mathcal{D}^{2ab+2m}_{2m, \psi^{\alpha}}(\CE_{\Delta(\tau, b)\otimes \sigma})$ be the $\psi^{\alpha}$-descent of $\CE_{\Delta(\tau, b)\otimes \sigma}$ (\cite[Section 3.2]{GRS11}). Then as a representation of $\wt{\Sp}_{2ab}(\BA)$, it is square-integrable and contains the whole space of the residual representation $\wt{\CE}_{\Delta(\tau,b)}$.

\textbf{Step (3).} $\wt{\CE}_{\Delta(\tau,b)}$ has a nonzero Fourier coefficient attached to the symplectic partition $[(a)^{2b}]$.

\textbf{Proof of Part (2) of Theorem \ref{main12}.}
From the results in \textbf{Steps (1)-(3)} above, we can see that $\CE_{\Delta(\tau, b+1)\otimes \sigma}$ has a nonzero Fourier coefficient attached to the composite partition $[(2m)1^{2ab}] \circ [(a)^{2b}]$ (for the definition of composite partitions and the attached Fourier coefficients, we refer to \cite[Section 1]{GRS03}).
Therefore, by \cite[Lemma 3.1]{JL15a} or \cite[Lemma 2.6]{GRS03},
$\CE_{\Delta(\tau, b)\otimes \sigma}$ has a nonzero Fourier coefficient attached to $[(2m)(a)^{2b}]$, which completes the proof of the Part (2) of Theorem \ref{main12}.


\subsection{Proof of Step (1)}

Note that by \cite[Lemma 1.1]{GRS03}, $\CE_{\Delta(\tau, b)\otimes \sigma}$ has a nonzero Fourier coefficient attached to the partition $[(2m)1^{2ab}]$ with respect to the character $\psi_{[(2m)1^{2ab}], \alpha}$ if and only if the $\psi^{\alpha}$-descent $\mathcal{D}^{2ab+2m}_{2m, \psi^{\alpha}}(\CE_{\Delta(\tau, b)\otimes \sigma})$ of $\CE_{\Delta(\tau, b)\otimes \sigma}$
is not identically zero as a representation of $\wt{\Sp}_{2ab}(\BA)$.

Recall that $P^{2l}_r = M^{2l}_r N^{2l}_r$ (with $1 \leq r \leq l$) is the standard parabolic subgroup
of ${\Sp}_{2l}$ with Levi part $M^{2l}_r$ isomorphic to $\GL_r \times \Sp_{2l-2r}$,
$N^{2l}_r$ is the unipotent radical.
$\wt{P}^{2l}_r(\BA)$ is the pre-image of $P^{2l}_r(\BA)=\wt{M}^{2l}_r(\BA) N^{2l}_r(\BA)$ in $\wt{\Sp}_{2l}(\BA)$.

Take any $\xi \in \CE_{\Delta(\tau, b)\otimes \sigma}$,
we will calculate the constant term of the Fourier-Jacobi coefficient $\CFJ^{\phi}_{\psi^{\alpha}_{m-1}}({\xi})$ along
${P}^{2ab}_r$, which is denoted by $\CC_{N^{2ab}_r}(\CFJ^{\phi}_{\psi^{\alpha}_{m-1}}({\xi}))$, where $1 \leq r \leq ab$.

By \cite[Theorem 7.8]{GRS11},
\begin{align}\label{sec4equ1}
\begin{split}
& \CC_{N^{2ab}_r}(\CFJ^{\phi}_{\psi^{\alpha}_{m-1}}({\xi}))\\
= \ & \sum_{k=0}^r \sum_{\gamma \in P^1_{r-k, 1^k}(F) \bs \GL_r(F)}
\int_{L(\BA)} \phi_1(i(\lambda)) \CFJ^{\phi_2}_{\psi^{\alpha}_{m-1+k}}
(\CC_{N^{2ab+2m}_{r-k}} ({\xi}))(\hat{\gamma} \lambda \beta) d \lambda.
\end{split}
\end{align}
We explain the notation used in \eqref{sec4equ1} as follows. $N^{2ab+2m}_{r-k}$ denotes the unipotent radical of the parabolic
subgroup $P^{2ab+2m}_{r-k}$ of ${\Sp}_{2ab+2m}$ with the Levi subgroup
$\GL_{r-k} \times \Sp_{2ab+2m-2r+2k}$.
$P^1_{r-k, 1^k}$ is a subgroup of $\GL_r$
consisting of matrices of the form
$\begin{pmatrix}
g & x\\
0 & z
\end{pmatrix}$, with $z \in U_k$, the standard maximal unipotent subgroup of $\GL_k$.
For $g \in \GL_j$, $j \leq ab+m$, $\hat{g}=\diag(g, I_{2ab+2m-2j}, g^*)$.
$L$ is a unipotent subgroup, consisting of matrices of the form
$\lambda = \begin{pmatrix}
I_r & 0\\
x & I_m
\end{pmatrix}^{\wedge}$, and $i(\lambda)$ is the last row of $x$.
$\beta=\begin{pmatrix}
0 & I_r\\
I_m & 0
\end{pmatrix}^{\wedge}$. We assume that
$\phi = \phi_1 \otimes \phi_2$, with
$\phi_1 \in \CS(\BA^r)$ and $\phi_2 \in \CS(\BA^{ab-r})$.
Finally, the Fourier-Jacobi coefficients have the following identity:
$$
\CFJ^{\phi_2}_{\psi^{\alpha}_{m-1+k}}
(\CC_{N^{2ab+2m}_{r-k}} ({\xi}))(\hat{\gamma} \lambda \beta):=
\CFJ^{\phi_2}_{\psi^{\alpha}_{m-1+k}}
(\CC_{N^{2ab+2m}_{r-k}} (\rho(\hat{\gamma} \lambda \beta){\xi}))(I),
$$
with $\rho(\hat{\gamma} \lambda \beta)$ denoting the right translation by $\hat{\gamma} \lambda \beta$, and the function is regarded
as taking first the constant term
$\CC_{N^{2ab+2m}_{r-k}} (\rho(\hat{\gamma} \lambda \beta){\xi})$, and then after restricted to ${\Sp}_{2ab+2m-2r+2k}(\BA)$,
taking the Fourier-Jacobi coefficient $\CFJ^{\phi_2}_{\psi^{\alpha}_{m-1+k}}$, which is a map
taking automorphic forms on ${\Sp}_{2ab+2m-2r+2k}(\BA)$
to those on $\wt{\Sp}_{2ab-2r}(\BA)$.

By the cuspidal support of ${\xi}$,
$\CC_{N^{2ab+2m}_{r-k}} ({\xi})$ is identically zero, unless $k=r$
or $r-k = la$, $1 \leq l \leq b$. When $k=r$, since
$[(2m+2r)1^{2ab-2r}]$ is bigger than $\eta_{\frak{so}_{2n+1}(\BC), \frak{sp}_{2n}(\BC)}(\ul{p}(\psi))$ under the lexicographical ordering, by \cite[Proposition 6.4]{JL15c} and \cite[Lemma 1.1]{GRS03}, $\CFJ^{\phi_2}_{\psi^{\alpha}_{m-1+r}}
({\xi})$ is identically zero, and hence
the corresponding term is zero.
When $r-k=la$, $1 \leq l \leq b$ and
$1 \leq k \leq r$,
then by Lemma \ref{constantterm}, after restricting to $\Sp_{2a(b-l)+2m}(\BA)$, $\CC_{N^{2ab+2m}_{r-k}} (\rho(\hat{\gamma} \lambda \beta){\xi})$ becomes a form in
$\CE_{\Delta(\tau,b-l) \otimes \sigma}$ whose Arthur parameter is $\psi' = (\tau, 2b-2l+1) \boxplus \boxplus_{i=2}^r (\tau_i,1)$.
Since
$[(2m+2k)1^{2a(b-l)-2k}]$ is bigger than $\eta_{\frak{so}_{2n'+1}(\BC), \frak{sp}_{2n'}(\BC)}(\ul{p}(\psi'))$ under the lexicographical ordering, where $2n'=2a(b-l)+2m$,
by \cite[Proposition 6.4]{JL15c} and \cite[Lemma 1.1]{GRS03}, it follows that  $\CFJ^{\phi_2}_{\psi^{\alpha}_{m-1+k}}
(\CC_{N^{2ab+2m}_{r-k}} (\rho(\hat{\gamma} \lambda \beta){\xi}))$ is also identically zero, and hence
the corresponding term is also zero.
Therefore, the only possibilities that $\CC_{N^{2ab}_r}(\CFJ^{\phi}_{\psi^{\alpha}_{m-1}}({\xi})) \neq 0$ are $r=la$, $1 \leq l \leq b$, and $k=0$. To prove that $\CFJ^{\phi}_{\psi^{\alpha}_{m-1}}({\xi})$ is not identically zero, we just have to show $\CC_{N^{2ab}_r}(\CFJ^{\phi}_{\psi^{\alpha}_{m-1}}({\xi})) \neq 0$ for some $r$.

Let $r=ab$, then
\begin{align}\label{sec4equ2}
\begin{split}
& \CC_{N^{2ab}_{ab}}(\CFJ^{\phi}_{\psi^{\alpha}_{m-1}}({\xi}))\\
= \ & \int_{L(\BA)} \phi_1(i(\lambda)) \CFJ^{\phi_2}_{\psi^{\alpha}_{m-1}}
(\CC_{N^{2ab+2m}_{ab}} ({\xi}))(\lambda \beta) d \lambda.
\end{split}
\end{align}
By Lemma \ref{constantterm},
when restricted to $\GL_{2ab}(\BA) \times \Sp_{2m}(\BA)$,
$$
\CC_{N^{2ab+2m}_{ab}} ({\xi}) \in \delta_{P^{2ab+2m}_{ab}}^{\frac{1}{2}}
\lvert \det \rvert^{-\frac{b+1}{2}} \Delta(\tau,b) \otimes \sigma.
$$
Clearly, the integral in \eqref{sec4equ2} is not identically zero if and only if $\sigma$ is $\psi^{\alpha}$-generic.
By assumption, $\sigma$ is $\psi^{\alpha}$-generic,  and hence $\CFJ^{\phi}_{\psi^{\alpha}_{m-1}}({\xi})$ is not identically zero.
Therefore, $\CE_{\Delta(\tau, b)\otimes \sigma}$ has a nonzero Fourier coefficient attached to the partition $[(2m)1^{2ab}]$ with respect to the character $\psi_{[(2m)1^{2ab}], \alpha}$.
This completes the proof of Step (1).


\subsection{Proof of Step (2)}

The proof of irreducibility of $\wt{\CE}_{\Delta(\tau, b)}$ is similar to that of $\wt{\CE}_{\Delta(\tau, 1)}$ which is given in the proof of Theorem 2.1 of \cite{GRS11}. To show the square-integrable residual representation $\wt{\CE}_{\Delta(\tau, b)}$ is irreducible,
it suffices to show that at each local place $v$,
\begin{equation} \label{equirr1}
\Ind_{\wt{P}_{ab}(F_v)}^{\wt{\Sp}_{2ab}(F_v)} \mu_{\psi_v^{-\alpha}} \Delta(\tau_v, b) \lvert \cdot \rvert^{\frac{b}{2}}
\end{equation}
has a unique irreducible quotient, where we assume that $\psi \cong \otimes_v \psi_v$,
$P_{ab}$ is the parabolic subgroup of ${\Sp}_{2ab}$ with Levi subgroup isomorphic to $\GL_{ab}$, and $\wt{P}_{ab}(F_v)$ is the pre-image of ${P}_{ab}(F_v)$ in $\wt{\Sp}_{2ab}(F_v)$.
Note that $\Delta(\tau_v, b)$ is the unique irreducible quotient of the following induced representation
\begin{equation*}
\Ind_{Q_{a^b}(F_v)}^{\GL_{ab}(F_v)}\tau_v \lvert \cdot \rvert^{\frac{b-1}{2}} \otimes \tau_v \lvert \cdot \rvert^{\frac{b-3}{2}} \otimes \cdots \otimes \tau_v \lvert \cdot \rvert^{\frac{1-b}{2}},
\end{equation*}
where $Q_{a^b}$ is the parabolic subgroup of $\GL_{ab}$ with Levi subgroup isomorphic to $\GL_{a}^{\times b}$.
Let $P_{a^b}$ be the parabolic subgroup of ${\Sp}_{2ab}$ with Levi subgroup isomorphic to $\GL_{a}^{\times b}$, and $\wt{P}_{a^b}(F_v)$ is the pre-image of ${P}_{a^b}(F_v)$ in $\wt{\Sp}_{2ab}(F_v)$.
We just have to show that the following induced representation has a unique irreducible quotient
\begin{equation}\label{equirr2}
\Ind_{\wt{P}_{a^b}(F_v)}^{\wt{\Sp}_{2ab}(F_v)} \mu_{\psi_v^{-\alpha}} \tau_v \lvert \cdot \rvert^{\frac{2b-1}{2}} \otimes \tau_v \lvert \cdot \rvert^{\frac{2b-3}{2}} \otimes \cdots \otimes \tau_v \lvert \cdot \rvert^{\frac{1}{2}}.
\end{equation}
Since $\tau_v$ is generic and unitary, by \cite{T86} and \cite{V86}, $\tau_v$ is fully parabolic induced from its Langlands data with exponents in the open interval $(-\frac{1}{2}, \frac{1}{2})$. Explicitly, we can assume that
$$\tau_v \cong \rho_1 \lvert \cdot \rvert^{\alpha_1} \times \rho_2 \lvert \cdot \rvert^{\alpha_2} \times \cdots \times \rho_r \lvert \cdot \rvert^{\alpha_r},$$
where $\rho_i$'s are tempered representations, $\alpha_i \in \BR$ and $\frac{1}{2} > \alpha_1 > \alpha_2 > \cdots > \alpha_r > -\frac{1}{2}$.
Therefore, the induced representation in \eqref{equirr2} can be written as follows:
\begin{align*}
& \mu_{\psi_v^{-\alpha}} \rho_1 \lvert \cdot \rvert^{\frac{2b-1}{2}+ \alpha_1} \times \rho_2 \lvert \cdot \rvert^{\frac{2b-1}{2}+\alpha_2} \times \cdots \times \rho_r \lvert \cdot \rvert^{\frac{2b-1}{2}+\alpha_r}\\
\times \ & \rho_1 \lvert \cdot \rvert^{\frac{2b-3}{2}+ \alpha_1} \times \rho_2 \lvert \cdot \rvert^{\frac{2b-3}{2}+\alpha_2} \times \cdots \times \rho_r \lvert \cdot \rvert^{\frac{2b-3}{2}+\alpha_r}\\
\times \ & \cdots\\
\times \ & \rho_1 \lvert \cdot \rvert^{\frac{1}{2}+ \alpha_1} \times \rho_2 \lvert \cdot \rvert^{\frac{1}{2}+\alpha_2} \times \cdots \times \rho_r \lvert \cdot \rvert^{\frac{1}{2}+\alpha_r} \rtimes 1_{\wt{Sp}_0(F_v)}.
\end{align*}
Since $\alpha_i \in \BR$ and $\frac{1}{2} > \alpha_1 > \alpha_2 > \cdots > \alpha_r > -\frac{1}{2}$, we can easily see that
the exponents satisfy
\begin{align*}
& \frac{2b-1}{2}+ \alpha_1 > \frac{2b-1}{2}+\alpha_2 > \cdots > \frac{2b-1}{2}+\alpha_r \\
> \ & \frac{2b-3}{2}+ \alpha_1 > \frac{2b-3}{2}+\alpha_2 > \cdots > \frac{2b-3}{2}+\alpha_r \\
> \ & \cdots \\
> \ & \frac{1}{2}+ \alpha_1 > \frac{1}{2}+\alpha_2 > \cdots > \frac{1}{2}+\alpha_r > 0.
\end{align*}
By Langlands classification of metaplectic groups (see \cite{BW00} and \cite{BJ13}), one can see that the induced representation in \eqref{equirr2} has a unique irreducible quotient which is the Langlands quotient.
This completes the proof of irreducibility of $\wt{\CE}_{\Delta(\tau, b)}$.

To prove the square-integrability of $\mathcal{D}^{2ab+2m}_{2m, \psi^{\alpha}}(\CE_{\Delta(\tau, b)\otimes \sigma})$,
we need to calculate
the automorphic exponent attached to the non-trivial constant term
considered in Step (1)
($r=ab$, for definition of automorphic exponent see \cite[\Rmnum{1}.3.3]{MW95}).
For this, we need to consider the action of
$$\ol{g}=\diag(g, g^*) \in \GL_{ab}(\BA) \times \wt{\Sp}_{0}(\BA).$$
Since $r=ab$, $\beta=\begin{pmatrix}
0 & I_{ab}\\
I_m & 0
\end{pmatrix}^{\wedge}$.
Let
$$\wt{g}:= \beta \diag(I_m, \ol{g}, I_m) \beta^{-1} = \diag(g, I_{2m}, g^*).$$
Then changing variables in
\eqref{sec4equ2} $\lambda \mapsto \wt{g} \lambda \wt{g}^{-1}$ will give
a Jacobian $\lvert \det(g) \rvert^{-m}$.
On the other hand, by \cite[Formula (1.4)]{GRS11},
the action of $\ol{g}$ on $\phi_1$ gives $\gamma_{\psi^{-\alpha}}(\det(g)) \lvert \det(g) \rvert^{\frac{1}{2}}$. Therefore, $\ol{g}$ acts by $\Delta(\tau,b)(g)$ with
character
\begin{align*}
& \delta_{P^{2ab+2m}_{ab}}^{\frac{1}{2}}(\wt{g})
\lvert \det(g) \rvert^{-\frac{b+1}{2}} \lvert \det(g) \rvert^{-m}
\gamma_{\psi^{-\alpha}}(\det(g)) \lvert \det(g) \rvert^{\frac{1}{2}}\\
= \ & \gamma_{\psi^{-\alpha}}(\det(g)) \delta_{P^{2ab}_{ab}}^{\frac{1}{2}}(\ol{g})
\lvert \det(g) \rvert^{-\frac{b}{2}}.
\end{align*}
Therefore, as a function on $\GL_{ab}(\BA) \times \wt{\Sp}_{0}(\BA)$,
\begin{align} \label{sec4equ3}
\begin{split}
& \CC_{N^{2ab}_{ab}}(\CFJ^{\phi}_{\psi^{\alpha}_{m-1}}({\xi}))\\
\in \ &
\gamma_{\psi^{-\alpha}} \delta_{P^{2ab}_{ab}}^{\frac{1}{2}}
\lvert \det(\cdot) \rvert^{-\frac{b}{2}} \Delta(\tau,b) \otimes 1_{\wt{\Sp}_{0}(\BA)}.
\end{split}
\end{align}
Since, the cuspidal exponent of $\Delta(\tau,b)$ is
$\{(\frac{1-b}{2}, \frac{3-b}{2}, \ldots, \frac{b-1}{2})\}$,
the cuspidal exponent of
$\CC_{N^{2ab}_{ab}}(\CFJ^{\phi}_{\psi^{\alpha}_{m-1}}({\xi}))$
is
$\{(\frac{1-2b}{2}, \frac{3-2b}{2}, \ldots, -\frac{1}{2})\}$.
Hence, by Langlands square-integrability criterion (\cite[Lemma \Rmnum{1}.4.11 ]{MW95}),
the automorphic
representation $\mathcal{D}^{2ab+2m}_{2m, \psi^{\alpha}}(\CE_{\Delta(\tau, b)\otimes \sigma})$ is square-integrable.

From \eqref{sec4equ3}, it is easy to see that as a representation of $\GL_{ab}(\BA) \times \wt{\Sp}_{0}(\BA)$,
\begin{align} \label{sec4equ4}
\CC_{N^{2ab}_{ab}}(\mathcal{D}^{2ab+2m}_{2m, \psi^{\alpha}}(\CE_{\Delta(\tau, b)\otimes \sigma}))
=
\gamma_{\psi^{-\alpha}} \delta_{P^{2ab}_{ab}}^{\frac{1}{2}}
\lvert \det(\cdot) \rvert^{-\frac{b}{2}} \Delta(\tau,b) \otimes 1_{\wt{\Sp}_{0}(\BA)}.
\end{align}
From the cuspidal support of the Speh residual representation $\Delta(\tau,b)$ of $\GL_{ab}(\BA)$, one can now easily see that
\begin{align*}
& \ \CC_{N^{2ab}_{a^b}}(\mathcal{D}^{2ab+2m}_{2m, \psi^{\alpha}}(\CE_{\Delta(\tau, b)\otimes \sigma}))\\
= & \ \gamma_{\psi^{-\alpha}} \delta_{P^{2ab}_{a^b}}^{\frac{1}{2}}
\tau \lvert \cdot \rvert^{\frac{1-2b}{2}} \otimes \tau \lvert \cdot \rvert^{\frac{3-2b}{2}} \otimes \cdots \otimes \tau \lvert \cdot \rvert^{-\frac{1}{2}} \otimes 1_{\wt{\Sp}_{0}(\BA)},
\end{align*}
where $N^{2ab}_{a^b}$ is the unipotent radical of the parabolic subgroup $P^{2ab}_{a^b}$ with Levi isomorphic to $\GL_a^{\times b}$. By \cite[Corollary 3.14 (ii)]{MW95}, any non-cuspidal irreducible summand of
$\mathcal{D}^{2ab+2m}_{2m, \psi^{\alpha}}(\CE_{\Delta(\tau, b)\otimes \sigma})$ must be contained in the space $\wt{\CE}_{\tau^{\otimes b}, \Lambda}$, which is the residual representation
generated by residues of the Eisenstein series associated to the induced representation
$$\Ind_{\wt{P}^{2ab}_{a^b}(\BA)}^{\wt{\Sp}_{2ab}(\BA)} \gamma_{\psi^{-\alpha}}
\tau \lvert \cdot \rvert^{s_1} \otimes \tau \lvert \cdot \rvert^{s_2} \otimes \cdots \otimes \tau \lvert \cdot \rvert^{s_b},$$
at the point $\Lambda=\{\frac{1-2b}{2}, \frac{3-2b}{2}, \cdots, \frac{-1}{2}\}$.
Since the Speh residual representation $\Delta(\tau,b)$ of $\GL_{ab}(\BA)$ is irreducible, by taking residues in stages, one can easily see that the space of the residual representation $\wt{\CE}_{\tau^{\otimes b}, \Lambda}$ is exactly identical to that of $\wt{\CE}_{\Delta(\tau,b)}$.
Therefore, any non-cuspidal irreducible summand of
$\mathcal{D}^{2ab+2m}_{2m, \psi^{\alpha}}(\CE_{\Delta(\tau, b)\otimes \sigma})$ must be contained in the space
$\wt{\CE}_{\Delta(\tau,b)}$. Hence, the descent representation $\mathcal{D}^{2ab+2m}_{2m, \psi^{\alpha}}(\CE_{\Delta(\tau, b)\otimes \sigma})$
has a non-trivial intersection with the space of the residual representation $\wt{\CE}_{\Delta(\tau,b)}$. Since we have seen that $\wt{\CE}_{\Delta(\tau,b)}$ is irreducible, $\mathcal{D}^{2ab+2m}_{2m, \psi^{\alpha}}(\CE_{\Delta(\tau, b)\otimes \sigma})$
must contain the whole space of the residual representation $\wt{\CE}_{\Delta(\tau,b)}$.
This completes the proof of Step (2).


\subsection{Proof of Step (3)}

The proof of that $\wt{\CE}_{\Delta(\tau,b)}$ has a nonzero Fourier coefficient attached to the symplectic partition $[(a)^{2b}]$,
is very similar to the proof of \cite[Theorem 4.2.2]{L13b}, if $a$ is even. The idea is to apply Lemma
\ref{nvequ2} repeatedly and use induction on $b$. Note that the case of $\wt{\CE}_{\Delta(\tau,1)}$ has already been proved in
\cite[Theorem 8.1]{GRS11}.
We omit the details here for this case.

In the following, we assume that $a=2k+1$ and prove $\wt{\CE}_{\Delta(\tau,b)}$ has a nonzero Fourier coefficient attached to the symplectic partition $\ul{p}:=[(2k+1)^{2b}]$ by induction on $b$.
When $b=1$, it has been proved in \cite[Theorem 8.2]{GRS11}, we will use similar idea here.
Assume that $\wt{\CE}_{\Delta(\tau,b-1)}$ has a nonzero Fourier coefficient attached to the symplectic partition $[(2k+1)^{2b-2}]$.

Take any $\varphi \in \wt{\CE}_{\Delta(\tau,b)}$, its Fourier coefficients
attached to $\ul{p}$ are of the following form
\begin{equation}\label{case1step3equ1}
\varphi^{\psi_{\ul{p}}}(g) = \int_{[V_{\ul{p},2}]}
\varphi(vg) \psi_{\ul{p}}^{-1}(v)dv.
\end{equation}
For definitions of the unipotent group $V_{\ul{p},2}$
and its character $\psi_{\ul{p}}$, see \cite[Section 2]{JL15a}.

Note that the one-dimensional torus $\CH_{\ul{p}}$ defined in \cite[(2.1)]{JL15a} has elements of the following form
$\CH_{\ul{p}}(t)=\diag(A(t), A(t), \ldots, A(t))$, where
$A(t)=\diag(t^{2k}, t^{2k-2}, \ldots, t^{-2k})$, and there are $2b$-copies of $A(t)$. Also note that the group $L_{\ul{p}}(\BA)$ defined in \cite[Section 2]{JL15a} is isomorphic to $\GL_{2b}^{2k+1}(\BA)$, and the stabilizer of the character $\psi_{\ul{p}}$ in $L_{\ul{p}}$ is isomorphic to the diagonal embedding $\wt{\Sp}_{2b}^{\Delta}(\BA)$. Let $\iota$ be this diagonal embedding. Let $N=\left\{n(x):=\begin{pmatrix}
1 & 0 & x\\
0 & I_{2b-2} & 0\\
0 & 0 & 1
\end{pmatrix}\right\}$, then
\begin{equation}\label{case1step3equ2}
\iota(N)=\left\{\iota(n(x))= \begin{pmatrix}
I_{2k+1} & 0 & x I_{2k+1}\\
0 & I_{(2k+1)(2b-2)} & 0\\
0 & 0 & I_{2k+1}
\end{pmatrix}\right\}.
\end{equation}

To show the integral in \eqref{case1step3equ1} is non-vanishing, it suffices to show that the following integral is non-vanishing:
\begin{equation}\label{case1step3equ3}
\int_{F \bs \BA} \int_{[V_{\ul{p},2}]}
\varphi(vn(x)g) \psi_{\ul{p}}^{-1}(v)dvdx.
\end{equation}

Let $\omega$ be a Weyl element which sends $\CH_{\ul{p}}(t)$ to the following torus element
\begin{equation*}
\diag(A(t), t^{2k}I_{2b-2}, t^{2k-2}I_{2b-2}, \ldots, t^{-2k}I_{2b-2}, A(t)).
\end{equation*}
Then $\omega$ has the form $\diag(I_{2k+1}, \omega_1, I_{2k+1})$.
Conjugating from left by $\omega$, the integral in \eqref{case1step3equ3} becomes
\begin{equation}\label{case1step3equ4}
\int_{[W]}
\varphi(w\omega g) \psi_{W}^{-1}(w)dw,
\end{equation}
where $W=\omega V_{\ul{p},2} \iota(N) \omega^{-1}$,
and $\psi_{W}(w) = \psi_{\ul{p}}(\omega^{-1} w \omega)$.
Then elements of $W$ have the following form:
\begin{equation}\label{case1step3equ5}
w =
\begin{pmatrix}
z_{2k+1} & q_1 & q_2\\
0 & w' & q_1^*\\
0 & 0 & z_{2k+1}^*
\end{pmatrix}
\begin{pmatrix}
I_{2k+1} & 0 & 0\\
p_1 & I_{(2b-2)(2k+1)} & 0\\
p_2 & p_1^* & I_{2k+1}
\end{pmatrix},
\end{equation}
where $z_{2k+1} \in V_{2k+1}$, the standard maximal unipotent subgroup of $\GL_{2k+1}$;
$w' \in \omega_1 V_{[(2k+1)^{2b-2}],2} \omega_1^{-1}$;
$q_1 \in M_{(2k+1)\times ((2b-2)(2k+1))}$
with certain conditions; $p_1 \in M_{((2b-2)(2m+1)) \times (2m+1)}$, with certain conditions;
$q_2 \in M_{(2k+1) \times (2k+1)}$, symmetric with respect to the secondary diagonal, such that
$q_2(i,j)=0$, for $1 \leq j < i \leq 2k+1$, and
$q_2(1,1)=q_2(2,2)=\cdots=q_2(2k+1,2k+1)$;
$p_2 \in M_{(2k+1) \times (2k+1)}$, symmetric with respect to the secondary diagonal, such that
$p_2(i,j)=0$, for $1 \leq j \leq i \leq 2k+1$.

Next, as in the proof of Section 4, we apply Lemma \ref{nvequ2} to fill the zero entries in $q_1, q_2$ using the nonzero entries in $p_1, p_2$.
Similarly, to proceed, we need to define a sequence of one-dimensional root subgroups and put them in a correct order.

For $1 \leq i \leq k$ and $1 \leq j \leq i $, let
$\alpha^i_j = e_i + e_{2k+1-i+j}$.
For $1 \leq i \leq k$ and $i+1\leq j \leq  i+(2b-2)i$, let $\alpha^i_j = e_i - e_{(2k+1)+(2b-2)i-(j-1)}$.
For $k+1 \leq i \leq 2k$ and $1 \leq j \leq 2k+1-i$, let
$\alpha^i_j = e_i + e_{i+j}$.
For $k+1 \leq i \leq 2k$ and $(2k+1-i)+1 \leq j \leq (2k+1-i) + ((b-1)+(2b-2)(i-k-1))$, let $\alpha^i_j = e_i + e_{(2k+1)b-(b-1)-(2b-2)(i-k-1)+j}$. Finally, for $k+1 \leq i \leq 2k$ and
$$
(2k+1-i) + ((b-1)+(2b-2)(i-k-1)) + 1 \leq j \leq  (2k+1-i) + (2b-2)i,
$$
let $\alpha^i_j = e_i - e_{(2k+1)b-(j-1)}$.
For the above roots $\alpha^i_j$, let $X_{\alpha^i_j}$ be the corresponding one-dimensional root subgroup.

For $1 \leq i \leq k$ and $1 \leq j \leq i $, let
$\beta^i_j = -e_{2k+1-i+j}-e_{i+1}$. For $1 \leq i \leq k$ and
$i+1\leq j \leq  i+(2b-2)i$, let $\beta^i_j = e_{(2k+1)+(2b-2)i-(j-1)}-e_{i+1}$.
For $k+1 \leq i \leq 2k$ and $1 \leq j \leq 2k+1-i$, let
$\beta^i_j = -e_{i+j}-e_{i+1}$. For $k+1 \leq i \leq 2k$ and $(2k+1-i)+1 \leq j \leq (2k+1-i) + ((b-1)+(2b-2)(i-k-1))$, let $\beta^i_j = - e_{(2k+1)b-(b-1)-(2b-2)(i-k-1)+j}-e_{i+1}$. Finally,   for $k+1 \leq i \leq 2k$ and
$$
(2k+1-i) + ((b-1)+(2b-2)(i-k-1)) + 1 \leq j \leq  (2k+1-i) + (2b-2)i,
$$
let $\beta^i_j = e_{(2k+1)b-(j-1)}-e_{i+1}$.
For the above roots $\beta^i_j$, let $X_{\beta^i_j}$ be the corresponding one-dimensional root subgroup.

For $1 \leq i \leq k$, let $m_i = i+(2b-2)i$; and for $k+1 \leq i \leq 2k$, let $m_{i} = (2k+1-i)+(2b-2)i$.
Let $\wt{W}$ be the subgroup of $W$ with elements of the form as in \eqref{case1step3equ5}, but with the $p_1$ and $p_2$ parts zero. Let $\psi_{\wt{W}} = \psi_W|_{\wt{W}}$. For any subgroup of $W$ containing $\wt{W}$, we automatically extend $\psi_{\wt{W}}$ trivially to this subgroup and still denote the character by $\psi_{\wt{W}}$.

Next, we will apply Lemma \ref{nvequ2} to a sequence of quadruples. For $i$ goes from 1 to $k+1$, one can see that the following quadruple satisfies all the conditions for Lemma \ref{nvequ2}:
\begin{align*}
(\wt{W}_i, \psi_{\wt{W}}, \{X_{\alpha^i_j}\}_{j=1}^{m_i}, \{X_{\beta^i_j}\}_{j=1}^{m_i}),
\end{align*}
where
$$\wt{W}_i = \prod_{s=1}^{i-1} \prod_{j=1}^{m_s} X_{\alpha^s_j} \wt{W} \prod_{l=i+1}^{2k} \prod_{j=1}^{m_l} X_{\beta^l_j}.$$
Applying Lemma \ref{nvequ2}, one can see that the integral in \eqref{case1step3equ4} is non-vanishing if and only if the following integral is non-vanishing:
\begin{equation}\label{case1step3equ6}
\int_{[\wt{W}_{i}']} \varphi(w \omega g) \psi^{-1}_{\wt{W}_{i}'}(w) dw,
\end{equation}
where
\begin{equation}\label{case1step3equ7}
\wt{W}_{i}' = \prod_{s=1}^{i} \prod_{j=1}^{m_s} X_{\alpha^s_j} \wt{W} \prod_{l=i+1}^{2k} \prod_{j=1}^{m_l} X_{\beta^l_j},
\end{equation}
and $\psi_{\wt{W}_{i}'}$ is extended from
$\psi_{\wt{W}}$ trivially.

For $i$ goes from $k+2$ to $2k$, before applying Lemma \ref{nvequ2} repeatedly to certain sequence of quadruples as above, we need to take Fourier expansion of the resulting integral at the end of the step $i-1$ along $X_{e_i + e_{i}}$
(at the end of step $k+1$, one get the integral in \eqref{case1step3equ6} with $i=k+1$ there, at the end of step $s$, $k+2 \leq s \leq 2k-1$, one would get the integral in \eqref{case1step3equ8}). Under the action of $\GL_1$, we get two kinds of Fourier coefficients corresponding to the two orbits of the dual of $[X_{e_i + e_{i}}]$: the trivial one and the
non-trivial one. It turns out that any Fourier coefficient corresponding to the non-trivial orbit contains an inner
integral which is exactly the Fourier coefficients attached to the partition $[(2i)1^{(2k+1)(2b)-2i}]$, which is identically zero by \cite[Proposition 6.4]{JL15c}, since $i \geq k+2$.
Therefore only the Fourier coefficient attached to the trivial orbit survives.

After taking Fourier expansion of the resulting integral at the end of the step $i-1$ along $X_{e_i + e_{i}}$ as above, one can see that the following quadruple satisfies all the conditions for Lemma \ref{nvequ2}:
\begin{align*}
(X_{e_i + e_{i}}\wt{W}_i, \psi_{\wt{W}}, \{X_{\alpha^i_j}\}_{j=1}^{m_i}, \{X_{\beta^i_1}\}_{j=1}^{m_i}),
\end{align*}
where
$$\wt{W}_i =\prod_{s=1}^{i-1} \prod_{j=1}^{m_s} X_{\alpha^s_j} \prod_{t=k+2}^{i-1} X_{e_t + e_{t}} \wt{W} \prod_{l=i+1}^{2k} \prod_{j=1}^{m_l} X_{\beta^l_j}.$$
Applying Lemma \ref{nvequ2}, we can see that the resulting integral at the end of the step $i-1$ is non-vanishing if and only if the following integral is non-vanishing:
\begin{equation}\label{case1step3equ8}
\int_{[\wt{W}_{i}']} \varphi(w \omega g) \psi^{-1}_{\wt{W}_{i}'}(w) dw,
\end{equation}
where
\begin{equation}\label{case1step3equ9}
\wt{W}_{i}' = \prod_{s=1}^{i} \prod_{j=1}^{m_s} X_{\alpha^s_j} \prod_{t=k+2}^{i} X_{e_t + e_{t}} \wt{W} \prod_{l=i+1}^{2k} \prod_{j=1}^{m_l} X_{\beta^l_j},
\end{equation}
and $\psi_{\wt{W}_{i}'}$ is the trivial extension of
$\psi_{\wt{W}}$.

One can see that elements of $\wt{W}_{2k}'$ have the following form:
\begin{equation}\label{case1step3equ10}
w=\begin{pmatrix}
z_{2k+1} & q_1 & q_2\\
0 & w' & q_1^*\\
0 & 0 & z_{2k+1}^*
\end{pmatrix},
\end{equation}
where $z_{2k+1} \in V_{2k+1}$, which is the standard maximal unipotent subgroup of $\GL_{2k+1}$; $w' \in \omega_1 V_{[(2k+1)^{2b-2}],2} \omega_1^{-1}$; $q_1 \in \rm{Mat}_{(2k+1) \times (2k+1)(2b-2)}$ with $q_1(2k+1,j)=0$ for $1 \leq j \leq (2k+1)(2b-2)$;
$q_2 \in \Mat_{(2k+1) \times (2k+1)}$, symmetric with respect to the secondary diagonal, with $q_2(2k+1,1)=0$. For $w\in {\wt{W}_{2k}'}$ of form in $\eqref{case1step3equ10}$,
$$\psi_{\wt{W}_{2k}'}(w)=\psi(\sum_{i=1}^{2k} z_{i,i+1})
\psi_{[(2k+1)^{2b-2}]}(\omega_1^{-1} w' \omega_1).$$

Now consider the Fourier expansion of the integral in \eqref{case1step3equ8} along the one-dimensional root subgroup $X_{2e_{2k+1}}$. By the same reason as above, only the Fourier coefficient corresponding to the trivial orbit of the dual of
$[X_{2e_{2k+1}}]$ survives, which is actually equal to the integral in \eqref{case1step3equ8} (with $i=2k$ there):
\begin{equation}\label{case1step3equ11}
\int_{[W_{2k+1}]} \varphi(w \omega g) \psi_{W_{2k+1}}(w)^{-1} dw,
\end{equation}
where elements in $W_{2k+1}$ have the same structure as in \eqref{case1step3equ10}, except that $q_2(2k+1,1)$ is not identically zero.

One can see that the integral in \eqref{case1step3equ11} has an inner integral which is exactly $\varphi^{\psi_{N_{1^{2k}}}}$, using notation in Lemma \ref{lem22} below. On the other hand, we know that by Lemma \ref{lem22} below,
$\varphi^{\psi_{N_{1^{2k}}}} = \varphi^{\wt{\psi}_{N_{1^{2k+1}}}}$.
Therefore, the
integral in \eqref{case1step3equ11}
becomes
\begin{equation}\label{case1step3equ12}
\int_{[W'_{2k+1}]} \varphi(w \omega g) \psi_{W'_{2k+1}}(w)^{-1} dw,
\end{equation}
where any element in $W'_{2k+1}$
has the following form:
$$
w = w(z_{2k+1}, w',q_1, q_2) =
\begin{pmatrix}
z_{2k+1} & q_1 & q_2\\
0 & w' & q_1^*\\
0 & 0 & z_{2k+1}^*
\end{pmatrix},
$$
where $z_{2k+1} \in V_{2k+1}$, the standard maximal unipotent subgroup of $\GL_{2k+1}$; $w' \in \omega_1 V_{[(2k+1)^{2b-2}],2} \omega_1^{-1}$; $q_1 \in \rm{Mat}_{(2k+1) \times (2k+1)(2b-2)}$;
$q_2 \in \Mat_{(2k+1) \times (2k+1)}$, symmetric with respect to the secondary diagonal. For $w\in {{W}_{2k+1}'}$ as above,
$$\psi_{{W}_{2k+1}'}(w)=\psi(\sum_{i=1}^{2k} z_{i,i+1})
\psi_{[(2k+1)^{2b-2}]}(\omega_1^{-1} w' \omega_1).$$

Hence, the integral in \eqref{case1step3equ12} can be written as
\begin{equation}\label{case1step3equ13}
\int_{W_{2k+1}''} \varphi_{{P}_{2k+1}}(w \omega g) \psi_{{W}_{2k+1}''}(w)^{-1} dw,
\end{equation}
where
$W_{2k+1}''$ is a subgroup of $W_{2k+1}'$ consisting of elements of the form
$w(z_{2k+1}, w',0, 0)$, $\psi_{W_{2k+1}''}=\psi_{W_{2k+1}'}|_{W_{2k+1}''}$,
and
$\varphi_{{P}_{2m+1}}$ is the constant term of
$\varphi$ along the parabolic subgroup $\wt{P}_{2k+1}(\BA)=\wt{M}_{2k+1}(\BA)N_{2k+1}(\BA)$ of $\wt{\Sp}_{2b(2k+1)}(\BA)$ with the
Levi subgroup isomorphic to $\GL_{2k+1}(\BA) \times \wt{\Sp}_{(2b-2)(2k+1)}(\BA)$.

By Lemma \ref{constantterm5} below,
$\varphi(w \omega g)_{\wt{P}_{2k+1}(\BA)}$ is an automorphic form
in $\gamma_{\psi^{-\alpha}}\tau \lvert \cdot \rvert^{-\frac{2b-1}{2}} \otimes \wt{\CE}_{\Delta(\tau,b-1)}$ when restricted to the Levi subgroup.
Note that the restriction of $\psi_{{W}_{2k+1}'}$
to the $z_{2k+1}$-part gives us a Whittaker coefficient of $\tau$, and the restriction
to the $w'$-part gives a Fourier coefficient of $\wt{\mathcal{E}}_{\Delta(\tau,b-1)}$
attached to the partition $[(2k+1)^{2b-2}]$ up to the conjugation of
the Weyl element $\omega_1$. On the other hand, $\tau$ is
generic, and by induction assumption, $\wt{\mathcal{E}}_{\Delta(\tau,b-1)}$
has a nonzero Fourier coefficient attached to the partition $[(2k+1)^{2b-2}]$.
Therefore, we can make the conclusion that
$\wt{\mathcal{E}}_{\Delta(\tau,b)}$
has a nonzero $\psi_{\underline{p}}$-Fourier coefficient
attached to the partition $[(2k+1)^{2b}]$.
This completes the proof of Step (3), up to Lemmas \ref{constantterm5}
and \ref{lem22}, which are stated below.

We remark that as Lemmas \ref{constantterm} and \ref{lem2}, Lemmas \ref{constantterm5}
and \ref{lem22} below are also analogues of \cite[Lemmas 4.2.4 and 4.2.6]{L13b}, with similar arguments, and hence we again only state them without proofs.

\begin{lem}\label{constantterm5}
Let $\wt{P}_{(2k+1)i}(\BA) = \wt{M}_{(2k+1)i}(\BA) N_{(2k+1)i}(\BA)$ with $1 \leq i \leq b$
be the parabolic subgroup of $\wt{\Sp}_{2b(2k+1)}(\BA)$ with the Levi part
$$
\wt{M}_{(2k+1)i}(\BA) \cong \GL_{(2k+1)i}(\BA) \times \wt{\Sp}_{(2k+1)(2b-2)}(\BA).
$$
Let $\varphi$
be an arbitrary automorphic form in $\wt{\CE}_{\Delta(\tau, b)}$.
Denote by $\varphi_{{P}_{(2k+1)i}}$ the constant term of $\varphi$
along ${P}_{(2k+1)i}$.
Then, for $1 \leq i \leq b$, $\varphi_{{P}_{(2k+1)i}}$ belongs to
\begin{align*}
\CA(N_{(2k+1)i}(\BA)\wt{M}_{(2k+1)i}(F)\bs \wt{\Sp}_{2b(2k+1)}(\BA))_{\gamma_{\psi^{-\alpha}}\Delta(\tau, i)\lvert \cdot \rvert^{-\frac{2b-i}{2}}
\otimes \wt{\mathcal{E}}_{\Delta(\tau,b-i)}}.
\end{align*}
\end{lem}


\begin{lem}\label{lem22}
Let $N_{1^{p}}(\BA)$ be the unipotent radical of the parabolic subgroup
$\wt{P}_{1^{p}}(\BA)$ of $\wt{\Sp}_{2b(2k+1)}(\BA)$ with Levi part isomorphic to
$\GL_1^{\times p}(\BA) \times \wt{\Sp}_{2b(2k+1)-2p}(\BA)$. Let
$$
\psi_{N_{1^{p}}}(n) := \psi(n_{1,2}+ \cdots + n_{p,p+1}),
$$
and
$$
\widetilde{\psi}_{N_{1^{p}}}(n) := \psi(n_{1,2}+ \cdots + n_{p-1,p}),
$$
be two characters of $N_{1^{p}}(\BA)$. For any automorphic form $\varphi \in \wt{\mathcal{E}}_{\Delta(\tau,b)}$,
define $\psi_{N_{1^{p}}}$ and $\widetilde{\psi}_{N_{1^{p}}}$-Fourier coefficients
as follows:
\begin{equation}\label{lem22equ1}
\varphi^{\psi_{N_{1^{p}}}}(g) := \int_{[N_{1^{p}}]} \varphi(ng)\psi_{N_{1^{p}}}(n)^{-1}
dn,
\end{equation}
\begin{equation}\label{lem22equ2}
\varphi^{\widetilde{\psi}_{N_{1^{p}}}}(g) := \int_{[N_{1^{p}}]} \varphi(ng)\widetilde{\psi}_{N_{1^{p}}}(n)^{-1}
du.
\end{equation}
Then $\varphi^{\psi_{N_{1^{p}}}} \equiv 0, \forall p \geq 2k+1$;
$\varphi^{\psi_{N_{1^{2k}}}} = \varphi^{\widetilde{\psi}_{N_{1^{2k+1}}}}$.
\end{lem}

\section{Proof of Theorem \ref{main22}}

In this section, we prove that $\CE_{\Delta(\tau, b+1)\otimes \sigma}$ has a nonzero Fourier coefficient attached to $[(2m+2a)(a)^{2b}]$. Assume that $\sigma$ is $\psi^{\alpha}$-generic with $\alpha \in F^*/(F^*)^2$.

Similarly as in the proof of Part (2) of Theorem \ref{main12} in Section 5, we separate the proof of Theorem \ref{main22} into following \textbf{two steps}:

\textbf{Step (1).} $\CE_{\Delta(\tau, b+1)\otimes \sigma}$ has a nonzero Fourier coefficient attached to the partition $[(2m+2a)1^{2ab}]$ with respect to the character $\psi_{[(2m+2a)1^{2ab}], \alpha}$ (for definition, see \cite[Section 2]{JL15a}).

\textbf{Step (2).}  Let $\mathcal{D}^{2a(b+1)+2m}_{2m, \psi^{\alpha}}(\CE_{\Delta(\tau, b+1)\otimes \sigma})$ be the $\psi^{\alpha}$-descent from representation $\CE_{\Delta(\tau, b+1)\otimes \sigma}$ of $\Sp_{2a(b+1)+2m}(\BA)$ to a representation of $\wt{\Sp}_{2ab}(\BA)$. Then it is square-integrable and contains the whole space of the residual representation $\wt{\CE}_{\Delta(\tau,b)}$ which is irreducible and constructed at the beginning of Section 5.

\textbf{Proof of Theorem \ref{main22}.}
First, recall from the Step (3) in the proof of Part (2) of Theorem \ref{main12} that $\wt{\CE}_{\Delta(\tau,b)}$ has a nonzero Fourier coefficient attached to the symplectic partition $[(a)^{2b}]$.
From the results in \textbf{Steps (1)-(2)} above, we can see that $\CE_{\Delta(\tau, b+1)\otimes \sigma}$ has a nonzero Fourier coefficient attached to the composite partition $[(2m+2a)1^{2ab}] \circ [(a)^{2b}]$ (for the definition of composite partitions and the attached Fourier coefficients, we refer to \cite[Section 1]{GRS03}).
Therefore, by \cite[Lemma 3.1]{JL15a} or \cite[Lemma 2.6]{GRS03}, $\CE_{\Delta(\tau, b+1)\otimes \sigma}$ has a nonzero Fourier coefficient attached to $[(2m+2a)(a)^{2b}]$, which completes the proof of Theorem \ref{main22}.

Before proving the above two steps, we record the following lemma which is analogous to Lemma \ref{constantterm},
whose proof will be omitted.

\begin{lem}\label{constantterm2}
Let $P_{ai} = M_{ai} N_{ai}$ with $1 \leq i \leq b+1$ be the parabolic subgroup of $\Sp_{2a(b+1)+2m}$ whose
Levi part
$M_{ai} \cong \GL_{ai} \times \Sp_{a(2b+2-2i)+2m}$.
Let $\varphi$
be an arbitrary automorphic form in $\CE_{\Delta(\tau, b+1)\otimes \sigma}$.
Denote by $\varphi_{P_{ai}}(g)$ the constant term of $\varphi$
along $P_{ai}$.
Then, for $1 \leq i \leq b+1$,
$$
\varphi_{P_{ai}} \in \CA(N_{ai}(\BA)M_{ai}(F)\bs \Sp_{2a(b+1)+2m}(\BA))_{\Delta(\tau, i)\lvert \cdot \rvert^{-\frac{2b+1-i}{2}}
\otimes \mathcal{E}_{\Delta(\tau,b+1-i) \otimes \sigma}}.
$$
Note that when $i=b$, $\mathcal{E}_{\Delta(\tau,b+1-i) \otimes \sigma} ={\CE}_{\tau \otimes \sigma}$, which is not a residual representation
as explained at the end of Section 2.2,
is nonzero and generic by \cite[Theorem 7.1.3]{Sh10}; and
when $i=b+1$, $\mathcal{E}_{\Delta(\tau,b+1-i) \otimes \sigma} = \sigma$.
\end{lem}

\subsection{Proof of Step (1)}

By \cite[Lemma 1.1]{GRS03}, $\CE_{\Delta(\tau, b+1)\otimes \sigma}$ has a nonzero Fourier coefficient attached to the partition $[(2m+2a)1^{2ab}]$ with respect to the character $\psi_{[(2m+2a)1^{2ab}], \alpha}$ if and only if the $\psi^{\alpha}$-descent $\mathcal{D}^{2a(b+1)+2m}_{2m+2a, \psi^{\alpha}}(\CE_{\Delta(\tau, b+1)\otimes \sigma})$ of $\CE_{\Delta(\tau, b+1)\otimes \sigma}$, which is a representation of $\wt{\Sp}_{2ab}(\BA)$,
is not identically zero.

Take any $\xi \in \CE_{\Delta(\tau, b+1)\otimes \sigma}$,
we will calculate the constant term of $\CFJ^{\phi}_{\psi^{\alpha}_{m+a-1}}({\xi})$ along the parabolic subgroup
$\wt{P}^{2ab}_r(\BA)=\wt{M}^{2ab}_r(\BA)N^{2ab}_r(\BA)$ of $\wt{\Sp}_{2ab}(\BA)$ with Levi isomorphic to $\GL_r(\BA) \times \wt{\Sp}_{2ab-2r}(\BA)$, $1 \leq r \leq ab$, which is denoted by $\CC_{N^{2ab}_r}(\CFJ^{\phi}_{\psi^{\alpha}_{m+a-1}}({\xi}))$.

By \cite[Theorem 7.8]{GRS11},
\begin{align}\label{main22equ1}
\begin{split}
& \CC_{N^{2ab}_r}(\CFJ^{\phi}_{\psi^{\alpha}_{m+a-1}}({\xi}))\\
= \ & \sum_{k=0}^r \sum_{\gamma \in P^1_{r-k, 1^k}(F) \bs \GL_r(F)}
\int_{L(\BA)} \phi_1(i(\lambda)) \CFJ^{\phi_2}_{\psi^{\alpha}_{m+a-1+k}}
(\CC_{N^{2a(b+1)+2m}_{r-k}} ({\xi}))(\hat{\gamma} \lambda \beta) d \lambda,
\end{split}
\end{align}
The notation in \eqref{main22equ1} are explained in order. $N^{2a(b+1)+2m}_{r-k}$ is the unipotent radical of the parabolic
subgroup $P^{2a(b+1)+2m}_{r-k}$ of ${\Sp}_{2a(b+1)+2m}$. $P^1_{r-k, 1^k}$ is a subgroup of $\GL_r$
consisting of matrices of the form
$\begin{pmatrix}
g & x\\
0 & z
\end{pmatrix}$, with $z \in U_k$, the standard maximal unipotent subgroup of $\GL_k$.
For $g \in \GL_j$, $j \leq a(b+1)+m$, $\hat{g}=\diag(g, I_{2a(b+1)+2m-2j}, g^*)$.
$L$ is a unipotent subgroup, consisting of matrices of the form
$\lambda = \begin{pmatrix}
I_r & 0\\
x & I_{m+a}
\end{pmatrix}^{\wedge}$, and $i(\lambda)$ is the last row of $x$, and
$\beta=\begin{pmatrix}
0 & I_r\\
I_{m+a} & 0
\end{pmatrix}^{\wedge}$. Finally the Schwartz function $\phi = \phi_1 \otimes \phi_2$ with
$\phi_1 \in \CS(\BA^r)$ and $\phi_2 \in \CS(\BA^{ab-r})$, and the function
$$
\CFJ^{\phi_2}_{\psi^{\alpha}_{m+a-1+k}}
(\CC_{N^{2a(b+1)+2m}_{r-k}} ({\xi}))(\hat{\gamma} \lambda \beta):=
\CFJ^{\phi_2}_{\psi^{\alpha}_{m+a-1+k}}
(\CC_{N^{2a(b+1)+2m}_{r-k}} (\rho(\hat{\gamma} \lambda \beta){\xi}))(I),$$
with $\rho(\hat{\gamma} \lambda \beta)$ denoting the right translation by $\hat{\gamma} \lambda \beta$, is a composition of
the restriction to ${\Sp}_{2a(b+1)+2m-2r+2k}(\BA)$ of
$\CC_{N^{2a(b+1)+2m}_{r-k}} (\rho(\hat{\gamma} \lambda \beta){\xi})$ with
the Fourier-Jacobi coefficient $\CFJ^{\phi_2}_{\psi^{\alpha}_{m+a-1+k}}$, which
takes automorphic forms on ${\Sp}_{2a(b+1)+2m-2r+2k}(\BA)$ to these on $\wt{\Sp}_{2ab-2r}(\BA)$.

By the cuspidal support of ${\xi}$,
$\CC_{N^{2a(b+1)+2m}_{r-k}} ({\xi})$ is identically zero, unless $k=r$
or $r-k = la$, with $1 \leq l \leq b+1$. When $k=r$, since
$[(2m+2a+2r)1^{2ab-2r}]$ is bigger than $\eta_{\frak{so}_{2n+1}(\BC), \frak{sp}_{2n}(\BC)}(\ul{p}(\psi))$ under the lexicographical ordering, by \cite[Proposition 6.4]{JL15c} and \cite[Lemma 1.1]{GRS03}, $\CFJ^{\phi_2}_{\psi^{\alpha}_{m+a-1+r}}
({\xi})$ is identically zero, hence
the corresponding term is zero.
When $r-k=la$, $1 \leq l \leq b+1$ and
$1 \leq k \leq r$,
by Lemma \ref{constantterm2}, after restricting to $\Sp_{2a(b+1-l)+2m}(\BA)$, $\CC_{N^{2a(b+1)+2m}_{r-k}} (\rho(\hat{\gamma} \lambda \beta){\xi})$ becomes a form in
$\CE_{\Delta(\tau,b+1-l) \otimes \sigma}$. Note that the Arthur parameter of $\CE_{\Delta(\tau,b+1-l) \otimes \sigma}$
is $\psi' = (\tau, 2b-2l+1) \boxplus (\tau,1)\boxplus_{i=3}^r (\tau_i,1)$ if $1 \leq l \leq b$, and
is $\psi' = \boxplus_{i=3}^r (\tau_i,1)$ if $l=b+1$.
Since
$[(2m+2k)1^{2a(b+1-l)-2k}]$ is bigger than $\eta_{\frak{so}_{2n'+1}(\BC), \frak{sp}_{2n'}(\BC)}(\ul{p}(\psi'))$ under the lexicographical ordering, where $2n'=2a(b+1-l)+2m$,
by \cite[Proposition 6.4]{JL15c} and \cite[Lemma 1.1]{GRS03}, $\CFJ^{\phi_2}_{\psi^{\alpha}_{m+a-1+k}}
(\CC_{N^{2a(b+1)+2m}_{r-k}} (\rho(\hat{\gamma} \lambda \beta){\xi}))$ is also identically zero. Hence
the corresponding term is also zero.

It follows that the only possibilities that $\CC_{N^{2a(b+1)+2m}_r}(\CFJ^{\phi}_{\psi^{\alpha}_{m+a-1}}({\xi})) \neq 0$ are
$r=la$ with $1 \leq l \leq b+1$, and $k=0$. To prove that $\CFJ^{\phi}_{\psi^{\alpha}_{m+a-1}}({\xi})$ is not identically zero,
we just have to show $\CC_{N^{2ab}_r}(\CFJ^{\phi}_{\psi^{\alpha}_{m+a-1}}({\xi})) \neq 0$ for some $r$.

Take $r=ab$. Then we have
\begin{align}\label{main22equ2}
\begin{split}
& \CC_{N^{2ab}_{ab}}(\CFJ^{\phi}_{\psi^{\alpha}_{m+a-1}}({\xi}))\\
= \ & \int_{L(\BA)} \phi_1(i(\lambda)) \CFJ^{\phi_2}_{\psi^{\alpha}_{m+a-1}}
(\CC_{N^{2a(b+1)+2m}_{ab}} ({\xi}))(\lambda \beta) d \lambda.
\end{split}
\end{align}
By Lemma \ref{constantterm2},
when restricted to $\GL_{2ab}(\BA) \times \Sp_{2m+2a}(\BA)$,
$$
\CC_{N^{2a(b+1)+2m}_{ab}} ({\xi}) \in \delta_{P^{2a(b+1)+2m}_{ab}}^{\frac{1}{2}}
\lvert \det \rvert^{-\frac{b+1}{2}} \Delta(\tau,b) \otimes (\CE_{\tau \otimes \sigma}),
$$
where $\CE_{\tau \otimes \sigma}$ is not a residual representation as explained at the end of Section 2.2.

Clearly, the integral in \eqref{main22equ2} is not identically zero if and only if $\CE_{\tau \otimes \sigma}$ is $\psi^{\alpha}$-generic.
Since by assumption, $\sigma$ is $\psi^{\alpha}$-generic, we have that $\CE_{\tau \otimes \sigma}$ is also $\psi^{\alpha}$-generic by
\cite[Theorem 7.1.3]{Sh10}. Hence, $\CFJ^{\phi}_{\psi^{\alpha}_{m+a-1}}({\xi})$ is not identically zero.
Therefore, $\CE_{\Delta(\tau, b)\otimes \sigma}$ has a nonzero Fourier coefficient attached to the partition $[(2m+2a)1^{2ab}]$
with respect to the character $\psi_{[(2m+2a)1^{2ab}], \alpha}$.
This completes the proof of Step (1).

\subsection{Proof of Step (2)}

To prove the square-integrability of the descent representation $\mathcal{D}^{2a(b+1)+2m}_{2m+2a, \psi^{\alpha}}(\CE_{\Delta(\tau, b+1)\otimes \sigma})$, as in Section 5.2,
we need to calculate
the automorphic exponent attached to the non-trivial constant term
considered in Step (1)
($r=ab$).
For this, we need to consider the action of
$$\ol{g}=\diag(g, g^*) \in \GL_{ab}(\BA) \times \wt{\Sp}_{0}(\BA).$$
Since $r=ab$, we have that $\beta=\begin{pmatrix}
0 & I_{ab}\\
I_{m+a} & 0
\end{pmatrix}^{\wedge}$.
Let
$$\wt{g}:=\beta \diag(I_{m+a}, \ol{g}, I_{m+a}) \beta^{-1} = \diag(g, I_{2m+2a}, g^*).$$
Then changing variables in
\eqref{sec4equ2} $\lambda \mapsto \wt{g} \lambda \wt{g}^{-1}$ will give
a Jacobian $\lvert \det(g) \rvert^{-m-a}$.
On the other hand, by \cite[Formula (1.4)]{GRS11},
the action of $\ol{g}$ on $\phi_1$ gives $\gamma_{\psi^{-\alpha}}(\det(g)) \lvert \det(g) \rvert^{\frac{1}{2}}$. Therefore, $\ol{g}$ acts by $\Delta(\tau,b)(g)$ with
character
\begin{align*}
& \delta_{P^{2a(b+1)+2m}_{ab}}^{\frac{1}{2}}(\wt{g})
\lvert \det(g) \rvert^{-\frac{b+1}{2}} \lvert \det(g) \rvert^{-m-a}
\gamma_{\psi^{-\alpha}}(\det(g)) \lvert \det(g) \rvert^{\frac{1}{2}}\\
= \ & \gamma_{\psi^{-\alpha}}(\det(g)) \delta_{P^{2ab}_{ab}}^{\frac{1}{2}}(\ol{g})
\lvert \det(g) \rvert^{-\frac{b}{2}}.
\end{align*}
Therefore, as a function on $\GL_{ab}(\BA) \times \wt{\Sp}_{0}(\BA)$,
\begin{align} \label{main22equ3}
\CC_{N^{2ab}_{ab}}(\CFJ^{\phi}_{\psi^{\alpha}_{m+a-1}}({\xi}))
\in
\gamma_{\psi^{-\alpha}} \delta_{P^{2ab}_{ab}}^{\frac{1}{2}}
\lvert \det(\cdot) \rvert^{-\frac{b}{2}} \Delta(\tau,b) \otimes 1_{\wt{\Sp}_{0}(\BA)}.
\end{align}
Since the cuspidal exponent of $\Delta(\tau,b)$ is
$\{(\frac{1-b}{2}, \frac{3-b}{2}, \ldots, \frac{b-1}{2})\}$,
the cuspidal exponent of
$\CC_{N^{2ab}_{ab}}(\CFJ^{\phi}_{\psi^{\alpha}_{m+a-1}}({\xi}))$
is
$\{(\frac{1-2b}{2}, \frac{3-2b}{2}, \ldots, -\frac{1}{2})\}$.
Hence, by Langlands square-integrability criterion (\cite[Lemma \Rmnum{1}.4.11]{MW95}),
the automorphic
representation $\mathcal{D}^{2a(b+1)+2m}_{2m+2a, \psi^{\alpha}}(\CE_{\Delta(\tau, b+1)\otimes \sigma})$ is square integrable.

From \eqref{main22equ3}, it follows that as a representation of $\GL_{ab}(\BA) \times \wt{\Sp}_{0}(\BA)$,
\begin{align} \label{main22equ4}
\begin{split}
& \CC_{N^{2ab}_{ab}}(\mathcal{D}^{2a(b+1)+2m}_{2m+2a, \psi^{\alpha}}(\CE_{\Delta(\tau, b+1)\otimes \sigma}))\\
= \ &
\gamma_{\psi^{-\alpha}} \delta_{P^{2ab}_{ab}}^{\frac{1}{2}}
\lvert \det(\cdot) \rvert^{-\frac{b}{2}} \Delta(\tau,b) \otimes 1_{\wt{\Sp}_{0}(\BA)}.
\end{split}
\end{align}
Therefore, a similar argument as in Section 5.2 implies that any non-cuspidal summand of
$\mathcal{D}^{2a(b+1)+2m}_{2m+2a, \psi^{\alpha}}(\CE_{\Delta(\tau, b+1)\otimes \sigma})$ must be an irreducible subrepresentation of
$\wt{\CE}_{\Delta(\tau,b)}$.
Hence, $\mathcal{D}^{2a(b+1)+2m}_{2m+2a, \psi^{\alpha}}(\CE_{\Delta(\tau, b+1)\otimes \sigma})$
has a non-trivial intersection with the space of the residual representation $\wt{\CE}_{\Delta(\tau,b)}$.
Since $\wt{\CE}_{\Delta(\tau,b)}$ is irreducible, $\mathcal{D}^{2a(b+1)+2m}_{2m+2a, \psi^{\alpha}}(\CE_{\Delta(\tau, b+1)\otimes \sigma})$
must contain the whole space of the residual representation $\wt{\CE}_{\Delta(\tau,b)}$.
This completes the proof of Step (2).

\section{Proof of Theorem \ref{main3}}

In this section, assuming that $a=2k$, $L(\frac{1}{2}, \tau \times \sigma) \neq 0$,
$\sigma \ncong 1_{\Sp_0(\BA)}$, and $\CE_{\tau \otimes \sigma}$ has a nonzero Fourier coefficient attached to the partition $[(2k+2m)(2k)]$, we prove that $\CE_{\Delta(\tau, b)\otimes \sigma}$ has a nonzero Fourier coefficient attached to the partition $[(2k+2m)(2k)^{2b-1}]$, for any $b \geq 1$.

Without loss of generality, by \cite[Lemma 3.1]{JL15a} or \cite[Lemma 2.6]{GRS03}, we may assume that
$\CE_{\tau \otimes \sigma}$ has a nonzero Fourier coefficient corresponding to the partition $[(2k+2m)1^{2k}]$ with respect to the character $\psi_{[(2k+2m)1^{2k}], \alpha}$, for some $\alpha \in F^*/(F^*)^2$. Then the $\psi^{\alpha}$-descent of $\CE_{\tau\otimes \sigma}$ is a generic representation of $\wt{\Sp}_{2k}(\BA)$.
Note that by the constant formula in \cite[Theorem 7.8]{GRS11}, one can easily see that this descent is also a cuspidal representation of $\wt{\Sp}_{2k}(\BA)$

Similarly as in previous sections, we separate the proof of Theorem \ref{main3} into following \textbf{three steps}:

\textbf{Step (1).} $\CE_{\Delta(\tau, b)\otimes \sigma}$ has a nonzero $\psi_{[(2k+2m)1^{2k(2b-1)}], \alpha}$-Fourier coefficient attached to the partition $[(2k+2m)1^{2k(2b-1)}]$ (for definition, see \cite[Section 2]{JL15a}).

\textbf{Step (2).}
Let $\wt{\sigma}$ be any irreducible subrepresentation of
the $\psi^{\alpha}$-descent of $\CE_{\tau\otimes \sigma}$. Then it is
a generic cuspidal representation of $\wt{\Sp}_{2k}(\BA)$ which is weakly lifting to $\tau$.
Using the theory of theta correspondence and the strong lifting from generic cuspidal representations of $\SO_{2n+1}(\BA)$ to automorphic representations of $\GL_{2n}(\BA)$, proved in \cite{JS03} (see also \cite{CKPSS04}), $\tau$ is also a strong lifting of $\wt{\sigma}$.

Define a residual representation $\wt{\CE}_{\Delta(\tau,b-1) \otimes \wt{\sigma}}$ as follows: for any
$$
\wt{\phi} \in \CA(N_{k(2b-1)}(\BA)\wt{M}_{k(2b-1)}(F) \bs \wt{\Sp}_{2k(2b-1)}(\BA))_{\gamma_{\psi^{-\alpha}} \Delta(\tau,b-1) \otimes \wt{\sigma}}
$$
one defines as in \cite{MW95}) the
residual Eisenstein series
$$
\wt{E}(\wt{\phi},s)(g)=\sum_{\gamma\in P_{k(2b-1)}(F)\bks \Sp_{2k(2b-1)}(F)}\lambda_s \wt{\phi}(\gamma g).
$$
It converges absolutely for real part of $s$ large and has meromorphic continuation to the whole complex plane $\BC$. By similar argument as that in \cite{JLZ13}, this Eisenstein series
has a simple pole at $\frac{b}{2}$, which is the right-most one.
Denote the representation generated by these residues at $s=\frac{b}{2}$
by $\wt{\CE}_{\Delta(\tau,b-1) \otimes \wt{\sigma}}$. This residual representation is square-integrable.
Since $\tau$ is also a strong lifting of $\wt{\sigma}$, the same argument as in Section 5.2 implies that $\wt{\CE}_{\Delta(\tau,b-1) \otimes \wt{\sigma}}$ is also irreducible (details will be omitted).

Let $\mathcal{D}^{4kb+2m}_{2k+2m, \psi^{\alpha}}(\CE_{\Delta(\tau, b)\otimes \sigma})$ be the $\psi^{\alpha}$-descent of $\CE_{\Delta(\tau, b)\otimes \sigma}$. Then as a representation of $\wt{\Sp}_{2k(2b-1)}(\BA)$, it is square-integrable and contains
the whole space of the residual representation $\wt{\CE}_{\Delta(\tau,b-1) \otimes \wt{\sigma}}$, where $\wt{\sigma}$ is an irreducible subrepresentation of
the $\psi^{\alpha}$-descent of $\CE_{\tau\otimes \sigma}$.

\textbf{Step (3).}
Let $\wt{\sigma}$ be any irreducible subrepresentation of
the $\psi^{\alpha}$-descent of $\CE_{\tau\otimes \sigma}$.
$\wt{\CE}_{\Delta(\tau,b-1) \otimes \wt{\sigma}}$ has a nonzero Fourier coefficient attached to the partition $[(2k)^{2b-1}]$.

\textbf{Proof of Theorem \ref{main3}.}
From the results in \textbf{Steps (1)-(3)} above, we can see that $\CE_{\Delta(\tau, b)\otimes \sigma}$ has a nonzero Fourier coefficient attached to the composite partition $[(2k+2m)1^{2k(2b-1)}] \circ [(2k)^{2b-1}]$ (for the definition of composite partitions and the attached Fourier coefficients, we refer to \cite[Section 1]{GRS03}).
Therefore, by \cite[Lemma 3.1]{JL15a} or \cite[Lemma 2.6]{GRS03},
$\CE_{\Delta(\tau, b)\otimes \sigma}$ has a nonzero Fourier coefficient attached to $[(2k+2m)(2k)^{2b-1}]$, which completes the proof of Theorem \ref{main3}.

Before proving the above three steps, we record the following lemma which is analogous to Lemmas \ref{constantterm} and \ref{constantterm2}.

\begin{lem}\label{constantterm3}
Let $P_{ai} = M_{ai} N_{ai}$ with $1 \leq i \leq b$ be the parabolic subgroup of $\Sp_{2ab+2m}$ with
Levi part $M_{ai} \cong \GL_{ai}\times \Sp_{a(2b-2i)+2m}$.
Let $\varphi$
be an arbitrary automorphic form in $\CE_{\Delta(\tau, b)\otimes \sigma}$.
Denote by $\varphi_{P_{ai}}(g)$ the constant term of $\varphi$
along $P_{ai}$.
Then, for $1 \leq i \leq b$,
$$
\varphi_{P_{ai}} \in \CA(N_{ai}(\BA)M_{ai}(F)\bs \Sp_{2ab+2m}(\BA))_{\Delta(\tau, i)\lvert \cdot \rvert^{-\frac{2b-i}{2}}
\otimes \mathcal{E}_{\Delta(\tau,b-i) \otimes \sigma}}.
$$
Note that when $i=b$, $\mathcal{E}_{\Delta(\tau,b-i) \otimes \sigma} = \sigma$.
\end{lem}

\subsection{Proof of Step (1)}

By \cite[Lemma 1.1]{GRS03}, $\CE_{\Delta(\tau, b)\otimes \sigma}$ has a nonzero Fourier coefficient attached to the partition $[(2k+2m)1^{2ab}]$ with respect to the character $\psi_{[(2k+2m)1^{2ab}], \alpha}$ if and only if the $\psi^{\alpha}$-descent $\mathcal{D}^{4kb+2m}_{2k+2m, \psi^{\alpha}}(\CE_{\Delta(\tau, b)\otimes \sigma})$ of $\CE_{\Delta(\tau, b)\otimes \sigma}$
is not identically zero, as a representation of $\wt{\Sp}_{2k(2b-1)}(\BA)$.

We calculate the constant term of $\CFJ^{\phi}_{\psi^{\alpha}_{k+m-1}}({\xi})$, for $\xi \in \CE_{\Delta(\tau, b)\otimes \sigma}$, along
the parabolic subgroup $\wt{P}^{2k(2b-1)}_r(\BA)=\wt{M}^{2k(2b-1)}_r(\BA)N^{2k(2b-1)}_r(\BA)$ of $\wt{\Sp}_{2k(2b-1)}(\BA)$ with Levi isomorphic to $\GL_r(\BA) \times \wt{\Sp}_{2k(2b-1)-2r}(\BA)$, which is denoted by $\CC_{N^{2k(2b-1)}_r}(\CFJ^{\phi}_{\psi^{\alpha}_{k+m-1}}({\xi}))$, where $1 \leq r \leq k(2b-1)$.

By \cite[Theorem 7.8]{GRS11},
\begin{align}\label{main3equ1}
\begin{split}
& \CC_{N^{2k(2b-1)}_r}(\CFJ^{\phi}_{\psi^{\alpha}_{k+m-1}}({\xi}))\\
= \ & \sum_{s=0}^r \sum_{\gamma \in P^1_{r-s, 1^s}(F) \bs \GL_r(F)}
\int_{L(\BA)} \phi_1(i(\lambda)) \CFJ^{\phi_2}_{\psi^{\alpha}_{k+m-1+s}}
(\CC_{N^{4kb+2m}_{r-s}} ({\xi}))(\hat{\gamma} \lambda \beta) d \lambda.
\end{split}
\end{align}
The notation in this formula are as follows. $N^{4kb+2m}_{r-s}$ is the unipotent radical of the parabolic
subgroup $P^{4kb+2m}_{r-s}$ of ${\Sp}_{4kb+2m}$ with Levi isomorphic to
$\GL_{r-s} \times \Sp_{4kb+2m-2r+2s}$. $P^1_{r-s, 1^s}$ is a subgroup of $\GL_r$
consisting of matrices of the form
$\begin{pmatrix}
g & x\\
0 & z
\end{pmatrix}$, with $z \in U_s$, the standard maximal unipotent subgroup of $\GL_s$.
For $g \in \GL_j$, $j \leq 2kb+m$, $\hat{g}=\diag(g, I_{4kb+2m-2j}, g^*)$.
$L$ is a unipotent subgroup, consisting of matrices of the form
$\lambda = \begin{pmatrix}
I_r & 0\\
x & I_{k+m}
\end{pmatrix}^{\wedge}$, and $i(\lambda)$ is the last row of $x$, and
$\beta=\begin{pmatrix}
0 & I_r\\
I_{k+m} & 0
\end{pmatrix}^{\wedge}$. The Schwartz function  $\phi = \phi_1 \otimes \phi_2$ with
$\phi_1 \in \CS(\BA^r)$ and $\phi_2 \in \CS(\BA^{k(2b-1)-r})$, and the function
$$
\CFJ^{\phi_2}_{\psi^{\alpha}_{k+m-1+s}}
(\CC_{N^{4kb+2m}_{r-s}} ({\xi}))(\hat{\gamma} \lambda \beta):=
\CFJ^{\phi_2}_{\psi^{\alpha}_{k+m-1+s}}
(\CC_{N^{4kb+2m}_{r-s}} (\rho(\hat{\gamma} \lambda \beta){\xi}))(I),$$
with $\rho(\hat{\gamma} \lambda \beta)$ denoting the right translation by $\hat{\gamma} \lambda \beta$, is a composition of
the restriction of $\CC_{N^{2ab+2m}_{r-k}} (\rho(\hat{\gamma} \lambda \beta){\xi})$ to ${\Sp}_{4kb+2m-2r+2s}(\BA)$ with
the Fourier-Jacobi coefficient $\CFJ^{\phi_2}_{\psi^{\alpha}_{k+m-1+s}}$, which
takes automorphic forms on ${\Sp}_{4kb+2m-2r+2s}(\BA)$
to those on $\wt{\Sp}_{4kb-2k-2r}(\BA)$.

By the cuspidal support of ${\xi}$,
$\CC_{N^{4kb+2m}_{r-s}} ({\xi})$ is identically zero, unless $s=r$
or $r-s = 2kl$, with $1 \leq l \leq b$. When $s=r$, since
$[(2k+2m+2r)1^{4kb-2k-2r}]$ is bigger than $\eta_{\frak{so}_{2n+1}(\BC), \frak{sp}_{2n}(\BC)}(\ul{p}(\psi))$ under the lexicographical ordering, by \cite[Proposition 6.4]{JL15c} and \cite[Lemma 1.1]{GRS03}, $\CFJ^{\phi_2}_{\psi^{\alpha}_{k+m-1+r}}
({\xi})$ is identically zero,  and hence
the corresponding term is zero.
When $r-s=la$, $1 \leq l \leq b$ and $1 \leq s \leq r$,
by Lemma \ref{constantterm3}, after restricting to $\Sp_{4k(b-l)+2m}(\BA)$, $\CC_{N^{4kb+2m}_{r-s}} (\rho(\hat{\gamma} \lambda \beta){\xi})$
becomes a form in
$\CE_{\Delta(\tau,b-l) \otimes \sigma}$. The Arthur parameter of $\CE_{\Delta(\tau,b-l) \otimes \sigma}$
is $\psi' = (\tau, 2b-2l) \boxplus \boxplus_{i=2}^r (\tau_i,1)$ .
Since
$[(2k+2m+2s)1^{4k(b-l)-2k-2s}]$ is bigger than $\eta_{\frak{so}_{2n'+1}(\BC), \frak{sp}_{2n'}(\BC)}(\ul{p}(\psi'))$ under the lexicographical ordering, where $2n'=4k(b-l)+2m$,
by \cite[Proposition 6.4]{JL15c} and \cite[Lemma 1.1]{GRS03},  $\CFJ^{\phi_2}_{\psi^{\alpha}_{k+m-1+s}}
(\CC_{N^{4kb+2m}_{r-s}} (\rho(\hat{\gamma} \lambda \beta){\xi}))$ is also identically zero, and hence
the corresponding term is also zero.
Therefore, the only possibilities that $\CC_{N^{2k(2b-1)}_r}(\CFJ^{\phi}_{\psi^{\alpha}_{k+m-1}}({\xi})) \neq 0$ are $r=2kl$, $1 \leq l \leq b$ and $s=0$. To prove that $\CFJ^{\phi}_{\psi^{\alpha}_{k+m-1}}({\xi})$ is not identically zero, we just have to show that  $\CC_{N^{2k(2b-1)}_r}(\CFJ^{\phi}_{\psi^{\alpha}_{k+m-1}}({\xi})) \neq 0$ for some $r$.

Taking $r=2k(b-1)$, we have
\begin{align}\label{main3equ2}
\begin{split}
& \CC_{N^{2k(2b-1)}_{2k(b-1)}}(\CFJ^{\phi}_{\psi^{\alpha}_{k+m-1}}({\xi}))\\
= \ & \int_{L(\BA)} \phi_1(i(\lambda)) \CFJ^{\phi_2}_{\psi^{\alpha}_{k+m-1}}
(\CC_{N^{4kb+2m}_{2k(b-1)}} ({\xi}))(\lambda \beta) d \lambda.
\end{split}
\end{align}
By Lemma \ref{constantterm3},
when restricted to $\GL_{2k(2b-2)}(\BA) \times \Sp_{4k+2m}(\BA)$,
$$
\CC_{N^{4kb+2m}_{2k(b-1)}} ({\xi}) \in \delta_{P^{4kb+2m}_{2k(b-1)}}^{\frac{1}{2}}
\lvert \det \rvert^{-\frac{b+1}{2}} \Delta(\tau,b-1) \otimes \CE_{\tau \otimes \sigma}.
$$
It follows that the integral in \eqref{main3equ2} is not identically zero if and only if $\CE_{\tau \otimes \sigma}$ has a nonzero
Fourier coefficient corresponding to the partition $[(2k+2m)1^{2k}]$ with respect to the character $\psi_{[(2k+2m)1^{2k}], \alpha}$.
Hence, by assumption, $\CFJ^{\phi}_{\psi^{\alpha}_{k+m-1}}({\xi})$ is not identically zero.
Therefore, $\CE_{\Delta(\tau, b)\otimes \sigma}$ has a nonzero Fourier coefficient attached to the partition $[(2k+2m)1^{2k(2b-1)}]$
with respect to the character $\psi_{[(2k+2m)1^{2k(2b-1)}], \alpha}$.
This completes the proof of Step (1).

\subsection{Proof of Step (2)}
To prove the square-integrability of the descent representation $\mathcal{D}^{4kb+2m}_{2k+2m, \psi^{\alpha}}(\CE_{\Delta(\tau, b)\otimes \sigma})$,
we need to calculate
the automorphic exponent attached to the non-trivial constant term
considered in Step (1) with $r=2k(b-1)$ (for definition of automorphic exponent see \cite[\Rmnum{1}.3.3]{MW95}).
For this, we need to consider the action of
$$\ol{g}=\diag(g, I_{2k}, g^*) \in \GL_{2k(b-1)}(\BA) \times \wt{\Sp}_{2k}(\BA).$$
Since $r=2k(b-1)$, $\beta=\begin{pmatrix}
0 & I_{2k(b-1)}\\
I_{k+m} & 0
\end{pmatrix}^{\wedge}$.
Let
$$\wt{g}:=\beta \diag(I_{k+m}, \ol{g}, I_{k+m}) \beta^{-1} = \diag(g, I_{4k+2m}, g^*).$$
Then changing variables in
\eqref{sec4equ2} $\lambda \mapsto \wt{g} \lambda \wt{g}^{-1}$ will give
a Jacobian $\lvert \det(g) \rvert^{-k-m}$.
On the other hand, by \cite[Formula (1.4)]{GRS11},
the action of $\ol{g}$ on $\phi_1$ gives $\gamma_{\psi^{-\alpha}}(\det(g)) \lvert \det(g) \rvert^{\frac{1}{2}}$. Therefore, $\ol{g}$ acts by $\Delta(\tau,b-1)(g)$ with
character
\begin{align*}
& \delta_{P^{4kb+2m}_{2k(b-1)}}^{\frac{1}{2}}
\lvert \det(g) \rvert^{-\frac{b+1}{2}} \lvert \det(g) \rvert^{-k-m}
\gamma_{\psi^{-\alpha}}(\det(g)) \lvert \det(g) \rvert^{\frac{1}{2}}\\
= \ & \gamma_{\psi^{-\alpha}}(\det(g)) \delta_{P^{2k(2b-1)}_{2k(b-1)}}^{\frac{1}{2}}(\ol{g})
\lvert \det(g) \rvert^{-\frac{b}{2}}.
\end{align*}
Therefore, combined with the calculation in Step (1), as a function on $\GL_{2k(b-1)}(\BA) \times \wt{\Sp}_{2k}(\BA)$,
\begin{align} \label{main3equ3}
\begin{split}
& \CC_{N^{2k(2b-1)}_{2k(b-1)}}(\CFJ^{\phi}_{\psi^{\alpha}_{k+m-1}}({\xi}))\\
\in \ &
\gamma_{\psi^{-\alpha}} \delta_{P^{2k(2b-1)}_{2k(b-1)}}^{\frac{1}{2}}
\lvert \det(\cdot) \rvert^{-\frac{b}{2}} \Delta(\tau,b-1) \otimes \mathcal{D}^{4k+2m}_{2k}(\CE_{\tau \otimes \sigma}).
\end{split}
\end{align}

Note that by the constant formula in \cite[Theorem 7.8]{GRS11}, one can easily see that $\mathcal{D}^{4k+2m}_{2k}(\CE_{\tau \otimes \sigma})$ is a cuspidal representation of $\wt{\Sp}_{2k}(\BA)$.
Since the cuspidal exponent of $\Delta(\tau,b-1)$ is
$\{(\frac{2-b}{2}, \frac{4-b}{2}, \ldots, \frac{b-2}{2})\}$,
the cuspidal exponent of
$\CC_{N^{2k(2b-1)}_{2k(b-1)}}(\CFJ^{\phi}_{\psi^{\alpha}_{k+m-1}}({\xi}))$
is
$\{(\frac{2-2b}{2}, \frac{4-2b}{2}, \ldots, -1)\}$.
Hence, by Langlands square-integrability criterion (\cite[Lemma \Rmnum{1}.4.11]{MW95}),
the automorphic
representation $\mathcal{D}^{4kb+2m}_{2k+2m, \psi^{\alpha}}(\CE_{\Delta(\tau, b)\otimes \sigma})$ is square-integrable.

From \eqref{main3equ3}, as a representation of $\GL_{2k(b-1)}(\BA) \times \wt{\Sp}_{2k}(\BA)$, we have
\begin{align} \label{main3equ4}
\begin{split}
& \CC_{N^{2k(2b-1)}_{2k(b-1)}}(\mathcal{D}^{4kb+2m}_{2k+2m, \psi^{\alpha}}(\CE_{\Delta(\tau, b)\otimes \sigma}))\\
= \ &
\gamma_{\psi^{-\alpha}} \delta_{P^{2k(2b-1)}_{2k(b-1)}}^{\frac{1}{2}}
\lvert \det(\cdot) \rvert^{-\frac{b}{2}} \Delta(\tau,b-1) \otimes \mathcal{D}^{4k+2m}_{2k, \psi^{\alpha}}(\CE_{\tau \otimes \sigma}).
\end{split}
\end{align}
Therefore, using a similar argument as in Section 5.2, one can see that $\mathcal{D}^{4kb+2m}_{2k+2m, \psi^{\alpha}}(\CE_{\Delta(\tau, b)\otimes \sigma})$
contains an irreducible subrepresentation of the residual representation $\wt{\CE}_{\Delta(\tau,b-1) \otimes \wt{\sigma}}$, where $\wt{\sigma}$ is an irreducible generic cuspidal representation of $\wt{\Sp}_{2k}(\BA)$ which is a subrepresentation of
the $\psi^{\alpha}$-descent of $\CE_{\tau\otimes \sigma}$, and is weakly lifting to $\tau$.
Since $\tau$ is also a strong lifting of $\wt{\sigma}$, a similar argument as in Section
5.2 implies that $\wt{\CE}_{\Delta(\tau,b-1) \otimes \wt{\sigma}}$ is irreducible.
Hence $\mathcal{D}^{4kb+2m}_{2k+2m, \psi^{\alpha}}(\CE_{\Delta(\tau, b)\otimes \sigma})$ must contain the whole space of residual representation $\wt{\CE}_{\Delta(\tau,b-1) \otimes \wt{\sigma}}$.
This completes the proof of Step (2).

\subsection{Proof of Step (3)}
Let $\wt{\sigma}$ be any irreducible subrepresentation of
the $\psi^{\alpha}$-descent of $\CE_{\tau\otimes \sigma}$, then it is a generic cuspidal representation of $\wt{\Sp}_{2k}(\BA)$.
Assume that $\wt{\sigma}$ is $\psi^{\beta}$-generic for some $\beta \in F^*/(F^*)^2$.

As in previous sections, we need to record the following lemma which is analogous to Lemma \ref{constantterm5}.

\begin{lem}\label{constantterm4}
Let $\wt{P}_{ai}(\BA) = \wt{M}_{ai}(\BA) N_{ai}(\BA)$ with $1 \leq i \leq b-1$ be the parabolic subgroup of $\wt{\Sp}_{2k(2b-1)}(\BA)$ with
Levi part
$$\wt{M}_{ai}(\BA) \cong \GL_{ai}(\BA) \times \wt{\Sp}_{2k(2b-1-2i)}(\BA).
$$
Let $\varphi$
be an arbitrary automorphic form in $\wt{\CE}_{\Delta(\tau,b-1) \otimes \wt{\sigma}}$.
Denote by $\varphi_{{P}_{ai}}(g)$ the constant term of $\varphi$
along ${P}_{ai}$.
Then, for $1 \leq i \leq b-1$,
$$
\varphi_{{P}_{ai}} \in \CA(N_{ai}(\BA)\wt{M}_{ai}(F)\bs \wt{\Sp}_{2k(2b-1)}(\BA))_{\gamma_{\psi^{-\alpha}}\Delta(\tau, i)\lvert \cdot \rvert^{-\frac{2b-1-i}{2}}
\otimes \wt{\mathcal{E}}_{\Delta(\tau,b-1-i) \otimes \wt{\sigma}}}.
$$
Note that when $i=b-1$, $\wt{\mathcal{E}}_{\Delta(\tau,b-1-i) \otimes \wt{\sigma}} = \wt{\sigma}$.
\end{lem}

First, we show that $\wt{\CE}_{\Delta(\tau,b-1) \otimes \wt{\sigma}}$ has a nonzero Fourier coefficient attached to the partition
$[(2k)1^{2k(2b-2)}]$ with respect to the character $\psi_{[(2k)1^{2k(2b-2)}], \beta}$.
By \cite[Lemma 1.1]{GRS03}, $\wt{\CE}_{\Delta(\tau,b-1) \otimes \wt{\sigma}}$ has a nonzero $\psi_{[(2k)1^{2k(2b-2)}], \beta}$-Fourier coefficient attached to the partition $[(2k)1^{2k(2b-2)}]$ if and only if the $\psi^{\beta}$-descent $\wt{\mathcal{D}}^{2k(2b-1)}_{2k, \psi^{\beta}}(\wt{\CE}_{\Delta(\tau,b-1) \otimes \wt{\sigma}})$ of $\wt{\CE}_{\Delta(\tau,b-1) \otimes \wt{\sigma}}$
is not identically zero, as a representation of ${\Sp}_{2k(2b-2)}(\BA)$.

Take any $\xi \in \wt{\CE}_{\Delta(\tau,b-1) \otimes \wt{\sigma}}$,
we will calculate the constant term of $\CFJ^{\phi}_{\psi^{\beta}_{k-1}}({\xi})$ along
${P}^{2k(2b-2)}_r$,  which is denoted by $\CC_{N^{2k(2b-2)}_r}(\CFJ^{\phi}_{\psi^{\beta}_{k-1}}({\xi}))$,
where $1 \leq r \leq k(2b-2)$. Recall that ${P}^{2k(2b-2)}_r={M}^{2k(2b-2)}_r{N}^{2k(2b-2)}_r$ is the parabolic subgroup of ${\Sp}_{2k(2b-2)}$ with Levi subgroup isomorphic to $\GL_r \times {\Sp}_{2k(2b-2)-2r}$.

By \cite[Theorem 7.8]{GRS11},
\begin{align}\label{main3equ5}
\begin{split}
& \CC_{N^{2k(2b-2)}_r}(\CFJ^{\phi}_{\psi^{\beta}_{k-1}}({\xi}))\\
= \ & \sum_{s=0}^r \sum_{\gamma \in P^1_{r-s, 1^s}(F) \bs \GL_r(F)}
\int_{L(\BA)} \phi_1(i(\lambda)) \CFJ^{\phi_2}_{\psi^{\beta}_{k-1+s}}
(\CC_{N^{2k(2b-1)}_{r-s}} ({\xi}))(\hat{\gamma} \lambda \eta) d \lambda.
\end{split}
\end{align}
Here are the notation in the formula. $N^{2k(2b-1)}_{r-s}(\BA)$ is the unipotent radical of the parabolic
subgroup $\wt{P}^{2k(2b-1)}_{r-s}(\BA)$ of $\wt{\Sp}_{2k(2b-1)}(\BA)$ with Levi subgroup isomorphic to
$\GL_{r-s}(\BA) \times \wt{\Sp}_{2k(2b-1)-2r+2s}(\BA)$. $P^1_{r-s, 1^s}$ is a subgroup of $\GL_r$
consisting of matrices of the form
$\begin{pmatrix}
g & x\\
0 & z
\end{pmatrix}$, with $z \in U_s$, the standard maximal unipotent subgroup of $\GL_s$.
For $g \in \GL_j$, $j \leq k(2b-1)$, $\hat{g}=\diag(g, I_{2k(2b-1)-2j}, g^*)$.
$L$ is a unipotent subgroup, consisting of matrices of the form
$\lambda = \begin{pmatrix}
I_r & 0\\
x & I_k
\end{pmatrix}^{\wedge}$, and $i(\lambda)$ is the last row of $x$, and
$\eta=\begin{pmatrix}
0 & I_r\\
I_k & 0
\end{pmatrix}^{\wedge}$. The Schwartz function $\phi = \phi_1 \otimes \phi_2$ with
$\phi_1 \in \CS(\BA^r)$ and $\phi_2 \in \CS(\BA^{k(2b-2)-r})$, and the function
$$\CFJ^{\phi_2}_{\psi^{\beta}_{k-1+s}}
(\CC_{N^{2k(2b-1)}_{r-s}} ({\xi}))(\hat{\gamma} \lambda \eta):=
\CFJ^{\phi_2}_{\psi^{\beta}_{k-1+s}}
(\CC_{N^{2k(2b-1)}_{r-s}} (\rho(\hat{\gamma} \lambda \eta){\xi}))(I),$$
with $\rho(\hat{\gamma} \lambda \eta)$ denoting the right translation by $\hat{\gamma} \lambda \eta$, is a composition of
the restriction to $\wt{\Sp}_{2k(2b-1)-2r+2s}(\BA)$ of
$\CC_{N^{2ab+2m}_{r-s}} (\rho(\hat{\gamma} \lambda \eta){\xi})$ with
the Fourier-Jacobi coefficient $\CFJ^{\phi_2}_{\psi^{\beta}_{k-1+s}}$, which
takes automorphic forms on $\wt{\Sp}_{2k(2b-1)-2r+2s}(\BA)$
to those on ${\Sp}_{2k(2b-2)-2r}(\BA)$.

By the cuspidal support of ${\xi}$,
$\CC_{N^{2k(2b-1)}_{r-s}} ({\xi})$ is identically zero, unless $s=r$
or $r-s = 2kl$, with $1 \leq l \leq b-1$. When $s=r$, from the structure of the unramified components of
the residual representation $\wt{\CE}_{\Delta(\tau,b-1) \otimes \wt{\sigma}}$,
by \cite[Lemma 3.2]{JL15c}, $\CFJ^{\phi_2}_{\psi^{\beta}_{k-1+r}}
({\xi})$ is identically zero, and hence
the corresponding term is zero.
When $r-s=2kl$, $1 \leq l \leq b-1$ and
$1 \leq s \leq r$,
then by Lemma \ref{constantterm4}, after restricting to $\wt{\Sp}_{2k(2b-1-2l)}(\BA)$, $\CC_{N^{2k(2b-1)}_{r-s}} (\rho(\hat{\gamma} \lambda \eta){\xi})$ becomes a form in
$\wt{\CE}_{\Delta(\tau,b-1-l) \otimes \wt{\sigma}}$.
From the structure of the unramified components of
the residual representation $\wt{\CE}_{\Delta(\tau,b-1-l) \otimes \wt{\sigma}}$,
by \cite[Lemma 3.2]{JL15c}, $FJ^{\phi_2}_{\psi^{\beta}_{k-1+s}}
(\CC_{N^{2k(2b-1)}_{r-s}} (\rho(\hat{\gamma} \lambda \eta){\xi}))$ is also identically zero, and hence
the corresponding term is also zero.
Therefore, the only possibilities that $\CC_{N^{2k(2b-2)}_r}(\CFJ^{\phi}_{\psi^{\beta}_{k-1}}({\xi})) \neq 0$ are $r=2kl$, $1 \leq l \leq b-1$, and $s=0$. To prove that $\CFJ^{\phi}_{\psi^{\beta}_{k-1}}({\xi})$ is not identically zero, we just have to show $\CC_{N^{2k(2b-2)}_r}(\CFJ^{\phi}_{\psi^{\beta}_{k-1}}({\xi})) \neq 0$ for some $r$.

Taking $r=2k(b-1)$, we have
\begin{align}\label{main3equ6}
\begin{split}
& \CC_{N^{2k(2b-2)}_r}(\CFJ^{\phi}_{\psi^{\beta}_{k-1}}({\xi}))\\
= \ & \int_{L(\BA)} \phi_1(i(\lambda)) \CFJ^{\phi_2}_{\psi^{\beta}_{k-1}}
(\CC_{N^{2k(2b-1)}_{2k(b-1)}} ({\xi}))(\lambda \eta) d \lambda.
\end{split}
\end{align}
By Lemma \ref{constantterm4},
when restricted to $\GL_{2k(b-1)}(\BA) \times \wt{\Sp}_{2k}(\BA)$,
$$
\CC_{N^{2k(2b-1)}_{2k(b-1)}} ({\xi}) \in \delta_{P^{2k(2b-1)}_{2k(b-1)}}^{\frac{1}{2}}
\lvert \det \rvert^{-\frac{b}{2}} \gamma_{\psi^{-\alpha}}\Delta(\tau,b-1) \otimes \wt{\sigma}.
$$
It is clear that the integral in \eqref{main3equ6} is not identically zero if and only if $\wt{\sigma}$ is $\psi^{\beta}$-generic. Hence, by assumption, $\CFJ^{\phi}_{\psi^{\alpha}_{k-1}}({\xi})$ is not identically zero.
Therefore, $\wt{\CE}_{\Delta(\tau, b-1)\otimes \wt{\sigma}}$ has a nonzero Fourier coefficient attached to the partition $[(2k)1^{2k(2b-2)}]$ with respect to the character $\psi_{[(2k)1^{2k(2b-2)}], \beta}$.

Next, we show that
the $\psi^{\beta}$-descent $\wt{\mathcal{D}}^{2k(2b-1)}_{2k, \psi^{\beta}}(\wt{\CE}_{\Delta(\tau, b-1)\otimes \wt{\sigma}})$ of $\wt{\CE}_{\Delta(\tau, b-1)\otimes \wt{\sigma}}$ is square-integrable and contains the whole space of the residual representation ${\CE}_{\Delta(\tau,b-1)}$ which is irreducible, as shown in \cite[Theorem 7.1]{L13a}.

To prove the square-integrability of $\wt{\mathcal{D}}^{2k(2b-1)}_{2k, \psi^{\beta}}(\wt{\CE}_{\Delta(\tau, b-1)\otimes \wt{\sigma}})$,
we need to calculate
the automorphic exponent attached to the non-trivial constant term
considered above
($r=2k(b-1)$).
For this, we need to consider the action of
$$\ol{g}=\diag(g, g^*) \in \GL_{2k(b-1)}(\BA) \times {\Sp}_{0}(\BA).$$
Since $r=2k(b-1)$, $\eta=\begin{pmatrix}
0 & I_{2k(b-1)}\\
I_k & 0
\end{pmatrix}^{\wedge}$.
Let
$$\wt{g}:=\eta \diag(I_k, \ol{g}, I_k) \eta^{-1} = \diag(g, I_{2k}, g^*).$$
Then changing variables in
\eqref{main3equ6} $\lambda \mapsto \wt{g} \lambda \wt{g}^{-1}$ will give
a Jacobian $\lvert \det(g) \rvert^{-k}$.
On the other hand, by \cite[Formula (1.4)]{GRS11},
the action of $\ol{g}$ on $\phi_1$ gives $\lvert \det(g) \rvert^{\frac{1}{2}}$. Therefore, $\ol{g}$ acts by $\Delta(\tau,b-1)(g)$ with
character
\begin{align*}
& \delta_{P^{2k(2b-1)}_{2k(b-1)}}^{\frac{1}{2}}
\lvert \det(g) \rvert^{-\frac{b}{2}} \lvert \det(g) \rvert^{-k}
\gamma_{\psi^{-\beta}}(\det(g)) \lvert \det(g) \rvert^{\frac{1}{2}}\\
= \ & \delta_{P^{2k(2b-2)}_{2k(b-1)}}^{\frac{1}{2}}(\ol{g})
\lvert \det(g) \rvert^{-\frac{b-1}{2}}.
\end{align*}
Therefore, as a function on $\GL_{2k(b-1)}(\BA) \times {\Sp}_{0}(\BA)$,
\begin{align} \label{main3equ7}
\CC_{N^{2k(2b-2)}_{2k(b-1)}}(\CFJ^{\phi}_{\psi^{\beta}_{k-1}}({\xi}))
\in
\delta_{P^{2k(2b-2)}_{2k(b-1)}}^{\frac{1}{2}}
\lvert \det(\cdot) \rvert^{-\frac{b-1}{2}} \Delta(\tau,b-1) \otimes 1_{{\Sp}_{0}(\BA)}.
\end{align}
Since the cuspidal exponent of $\Delta(\tau,b-1)$ is
$\{(\frac{2-b}{2}, \frac{4-b}{2}, \ldots, \frac{b-2}{2})\}$,
the cuspidal exponent of
$\CC_{N^{2k(2b-2)}_{2k(b-1)}}(\CFJ^{\phi}_{\psi^{\beta}_{k-1}}({\xi}))$
is
$\{(\frac{3-2b}{2}, \frac{5-2b}{2}, \ldots, -\frac{1}{2})\}$.
By Langlands square-integrability criterion (\cite[Lemma \Rmnum{1}.4.11]{MW95}),
the automorphic
representation $\wt{\mathcal{D}}^{2k(2b-1)}_{2k, \psi^{\beta}}(\wt{\CE}_{\Delta(\tau, b-1)\otimes \wt{\sigma}})$ is square integrable.

From \eqref{main3equ7}, it is easy to see that as a representation of $\GL_{2k(b-1)}(\BA) \times {\Sp}_{0}(\BA)$,
\begin{align} \label{main3equ8}
\begin{split}
& \CC_{N^{2k(2b-2)}_{2k(b-1)}}(\wt{\mathcal{D}}^{2k(2b-1)}_{2k, \psi^{\beta}}(\wt{\CE}_{\Delta(\tau, b-1)\otimes \wt{\sigma}}))\\
= \ &
\delta_{P^{2k(2b-2)}_{2k(b-1)}}^{\frac{1}{2}}
\lvert \det(\cdot) \rvert^{-\frac{b-1}{2}} \Delta(\tau,b-1) \otimes 1_{{\Sp}_{0}(\BA)}.
\end{split}
\end{align}
It follows that  $\wt{\mathcal{D}}^{2k(2b-1)}_{2k, \psi^{\beta}}(\wt{\CE}_{\Delta(\tau, b-1)\otimes \wt{\sigma}})$
has a non-trivial intersection with the space of the residual representation ${\CE}_{\Delta(\tau,b-1)}$. Since by \cite[Theorem 7.1, Part (2)]{L13a}, ${\CE}_{\Delta(\tau,b-1)}$ is irreducible, $\wt{\mathcal{D}}^{2k(2b-1)}_{2k, \psi^{\beta}}(\wt{\CE}_{\Delta(\tau, b-1)\otimes \wt{\sigma}})$
must contains the whole space of the residual representation ${\CE}_{\Delta(\tau,b-1)}$.
By \cite[Theorem 7.1, Part (3)]{L13a},
the descent $\wt{\mathcal{D}}^{2k(2b-1)}_{2k, \psi^{\beta}}(\wt{\CE}_{\Delta(\tau, b-1)\otimes \wt{\sigma}})$ is actually irreducible and equals to the residual representation ${\CE}_{\Delta(\tau,b-1)}$ identically.

By \cite[Theorem 4.2.2]{L13b}, we know that $\mathfrak{p}^m({\CE}_{\Delta(\tau,b-1)}) = \{[(2k)^{2b-2}]\}$. Therefore, by \cite[Lemma 3.1]{JL15a} or \cite[Lemma 2.6]{GRS03}, $\wt{\CE}_{\Delta(\tau, b-1)\otimes \wt{\sigma}}$ has a nonzero Fourier coefficient attached to the partition $[(2k)^{2b-1}]$.
This completes the proof of Step (3).

\end{document}